\documentclass[11pt]{amsart}
\usepackage{amssymb,mathrsfs}
\usepackage{indentfirst,latexsym,enumerate}
\usepackage{colortbl}
\usepackage{epsfig,subfigure,graphicx}
\usepackage{amsfonts}
\usepackage{amsmath}
\usepackage{amscd,amsthm}
\usepackage{enumitem}
\usepackage{needspace}
\usepackage{placeins}
\usepackage{changes}
\usepackage{bigints}
\usepackage[all]{xy}
\usepackage[noadjust]{cite}
\usepackage{amsaddr}
\usepackage{booktabs,array}
\usepackage[hidelinks]{hyperref}
\usepackage[norefs,nocites,ignoreunlbld]{refcheck}
\numberwithin{equation}{section}
\allowdisplaybreaks

\title[Planar central configurations of four equal masses]{Planar Central Configurations of Four Equal Masses:\ Classification and Morse Indices}
\author[X. Huangfu and H. Shen]{Xiaoxiang Huangfu, Huatong Shen}
\address{
  \vspace*{-0.1cm}
  \upshape
  \begin{tabular}{@{}l}
       School of Mathematics, Shandong University, Jinan 250100, China.
  \end{tabular}
}
\date{}

\theoremstyle{plain}
\newtheorem{theorem}{Theorem}[section]
\newtheorem{lemma}[theorem]{Lemma}
\newtheorem{corollary}[theorem]{Corollary}
\newtheorem{proposition}[theorem]{Proposition}
\theoremstyle{definition}

\theoremstyle{remark}

\newtheorem*{remark*}{Remark}

\begin{document}
\begin{abstract}
We study the planar central configurations of four equal positive masses under
homogeneous interactions with force exponent $A>2$. For $2<A\le 3$, we prove
that there are exactly four types of configurations up to similarity and
relabeling: the square, the centered equilateral triangle, a unique nontrivial
concave kite, and a unique collinear configuration. Combined with known results
for $A>3$, this completes the classification for all $A>2$ and gives an
affirmative answer to Problem~11 of Albouy, Cabral, and Santos~\cite{AlbouyCabralSantos2012}. We also
establish the nondegeneracy modulo rotation of these configurations and
determine their Morse indices, confirming Hampton's Conjecture~2~\cite{Hampton2019}. Our proofs are
based on analytic arguments and exact algebraic calculations.
\end{abstract}

\maketitle
\section{Introduction}
\label{sec:introduction}

Central configurations play an important role in the classical $N$-body problem, as they give rise to homographic solutions and, within the planar Newtonian setting, relative equilibria. Their existence, count, and geometric properties have been extensively investigated; see, for example, \cite{moeckel2015central}.

For the three-body problem, the central configurations are well understood and consist of the Euler collinear and Lagrange equilateral configurations. For four or more bodies, the classification problem becomes considerably more difficult. Hampton and Moeckel \cite{HamptonMoeckel2006} proved the finiteness of central configurations in the Newtonian four-body problem for arbitrary positive masses, up to symmetry. However, finiteness alone does not provide a complete classification. The collinear case goes back to Moulton~\cite{moulton1910straight}, who established the existence and uniqueness of collinear central configurations for each ordering of the masses.

In this paper, we consider planar central configurations of four equal positive
masses. Let $q_i\in\mathbb R^2$ denote the position of the $i$-th body, let
$r_{ij}=|q_i-q_j|$, and let $c$ be the center of mass. We study the equations
\[
 \lambda(q_i-c)
 =\sum_{j\ne i}\frac{q_i-q_j}{r_{ij}^{A}},
 \qquad A>2,
\]
where $\lambda>0$ is a constant. The case $A=3$ corresponds to the classical
Newtonian problem. Varying $A$ changes how the interaction between two bodies
depends on their distance. Since translations, rotations, and scalings do not
change the essential shape of a central configuration, we consider
configurations up to similarity and relabeling.

For the Newtonian four-body problem with equal masses,
Albouy~\cite{Albouy1995,Albouy1996} proved that every noncollinear central
configuration has a reflection symmetry and obtained a complete classification.
Up to similarity and relabeling, there are exactly four types: the square, the
equilateral triangle with one body at its center, a nontrivial concave kite,
and the collinear configuration. Here, a concave kite is a configuration in
which one body lies inside the triangle formed by the other three bodies,
and the configuration has a reflection axis.

The geometry and uniqueness of four-body central configurations have also been investigated under more general mass assumptions. Albouy, Fu, and Sun
\cite{albouy2008symmetry} studied the relationship between masses and reflection symmetry in convex configurations. More recently, Sun, Xie, and You
\cite{sun2023uniqueness} obtained computer-assisted uniqueness results for convex central configurations in specified regions of the mass parameter space.
Roberts~\cite{Roberts2025} gave a new proof of the uniqueness of convex kite configurations and studied the corresponding concave configurations. Related
symmetric configurations and generalized Laura--Andoyer equations have also been considered by Dias and Tsai~\cite{DiasTsai2026}.

Other approaches to the four-body central-configuration problem are based on mutual distances and angular coordinates. Albouy and Chenciner \cite{albouy1997probleme} developed a formulation of the $N$-body problem in terms of mutual distances. Cors and Roberts \cite{cors2012four} studied co-circular central configurations in the Newtonian four-body problem, including kite and isosceles-trapezoid configurations. For power-law potentials, Cors, Hall, and Roberts \cite{cors2014uniqueness} established uniqueness for co-circular central configurations whose center of mass coincides with the center of the circle. Santoprete \cite{santoprete2019planarity} derived alternative equations
for four-body central configurations using angles as coordinates.

For four equal masses, the nontrivial concave kite is distinguished by having exactly one reflection axis. A natural question is whether this uniqueness
persists when the exponent of the homogeneous interaction varies. Albouy, Cabral, and Santos \cite{AlbouyCabralSantos2012} raised this question
for exponents in a neighborhood of the Newtonian case. In their convention, the exponent is denoted by $a$, with $A=-2a$, so that the Newtonian value
corresponds to $a=-3/2$.

Hampton~\cite{Hampton2019} subsequently studied central configurations and their Morse indices for homogeneous interactions with varying exponents.
Morse theory has long been used in the study of central configurations; see, for example, Pacella~\cite{pacella1987central}. The Morse index counts the negative directions of the Hessian at a critical point, after the appropriate symmetry reduction.

For four equal masses, Hampton established the classification for $A\geq3$ and obtained the Morse polynomial
\[
M_A(s)=6+24s+20s^2.
\]
He also suggested that the same Morse polynomial should hold for every $A>2$. His method excludes bifurcations for $A>3$, but becomes difficult to apply in the interval $2<A<3$ because of the degeneracy at the limiting exponent $A=2$.

The main purpose of this paper is to complete the classification for $2<A<3$ and determine the corresponding Morse indices. Our first result is the following.

\begin{theorem}
\label{thm:planar-classification}
For every $2<A\le3$, there are exactly four types of planar central
configurations of four equal positive masses, up to similarity and relabeling:
the square, the equilateral triangle with the fourth mass at its center, a
unique nontrivial concave kite, and the unique collinear central configuration.
The nontrivial concave kite has exactly one reflection axis.
\end{theorem}

Combining Theorem~\ref{thm:planar-classification} with Hampton's result for
$A\ge3$, we obtain the classification for all $A>2$. In particular, we have the following uniqueness result.

\begin{corollary}
\label{cor-Newonian-case}
For every $a<-1$, where $A=-2a$, there exists exactly one planar central
configuration of four equal masses with precisely one reflection axis, up to
similarity and relabeling. In particular, this uniqueness holds on an interval
containing the Newtonian exponent $a=-3/2$.
\end{corollary}

We also study the nondegeneracy and Morse indices of these configurations.
Our second main result is as follows.

\begin{theorem}\label{Morse Theorem.}
    For every $2<A\le3$, all critical points of Hampton's four-equal-mass
function are nondegenerate modulo rotation, and their Morse polynomial is
\[
 M_A(s)=6+24s+20s^2.
\]
\end{theorem}

Together with Hampton's result for $A\ge3$, this establishes the same Morse polynomial for all $A>2$.

We briefly describe the main idea of the proof. By Albouy's symmetry result,
it is enough to consider a few symmetric types of configurations. The main
difficulty is to prove the uniqueness of the nontrivial concave kite. We first
reduce the central-configuration equations to two scalar equations and then
use one of them to obtain an equation in a single variable. At $A=2$, this
reduced equation admits an exact factorization. For $2<A\le3$, we prove that
the corresponding function is strictly convex on the relevant interval, which
gives the desired uniqueness. The remaining geometric cases are handled by
direct analytic arguments. Finally, we study the Hessian at each configuration
to determine its Morse index. Some of the inequalities in the proof are
verified using exact rational polynomial calculations, without interval
arithmetic.

The paper is organized as follows. Section~\ref{sec:reduction} derives the
reduced equations, and Section~\ref{sec:endpoint} studies the limiting case
$A=2$. Section~\ref{sec:concave} proves the uniqueness of the nontrivial concave
kite, while Section~\ref{sec:global-classification} completes the classification.
Section~\ref{sec:morse} determines the Morse indices, and
Appendix~\ref{app:exact-algebra} contains the exact algebraic calculations used
in the proof.

\section{Symmetry reduction and the reduced equations}
\label{sec:reduction}

In this section, we derive the reduced equations for concave central configurations of four equal masses. We first recall the central-configuration equations for homogeneous potentials and the symmetry result used in the reduction.

Let $q=(q_1,q_2,q_3,q_4)$ be a collision-free configuration, where $q_i\in\mathbb R^2$. We normalize the four masses to one and denote by $r_{ij}=|q_i-q_j|$ the mutual distances and by $c=\frac14\sum_{i=1}^4q_i$ the center of mass. Following Hampton's exponent convention, we consider $A>2$ and set $\gamma=A/2>1$, where $A=3$ corresponds to the Newtonian case.

The homogeneous potential is defined by
\[
U_A(q)=\frac{1}{A-2}
\sum_{1\leq i<j\leq4}r_{ij}^{2-A}.
\]
A central configuration satisfies
\begin{equation}
\lambda(q_i-c)
=\sum_{j\ne i}\frac{q_i-q_j}{r_{ij}^{A}},
\qquad i=1,\ldots,4,
\label{eq:cc}
\end{equation}
for some $\lambda>0$. Multiplying these equations by $q_i-c$ and summing over $i$, we obtain
\[
\lambda\sum_{i=1}^4|q_i-c|^2
=\sum_{1\leq i<j\leq4}r_{ij}^{2-A}.
\]
In particular, the multiplier $\lambda$ is positive. Moreover, translations, rotations and positive scalings preserve central configurations. Under the scaling $q_i\mapsto rq_i$, the multiplier changes to $r^{-A}\lambda$. We shall therefore study central configurations up to similarity and relabeling.

To apply the symmetry result of Albouy, we write the pair potential as a function of the squared distance $s=r^2$:
\[
\Phi_A(s)=\frac{s^{1-\gamma}}{A-2}.
\]
Direct differentiation gives
\[
\Phi_A''(s)=\frac{\gamma}{2}s^{-\gamma-1}>0,
\qquad
\Phi_A'''(s)
=-\frac{\gamma(\gamma+1)}2s^{-\gamma-2}<0.
\]
Thus the symmetry result of Albouy \cite[Propositions~1 and~2]{Albouy1995}
applies to the present homogeneous potential (see also \cite[Section~3.2]{Hampton2019}). In particular, every noncollinear planar central
configuration of four equal masses has at least one axis of reflection symmetry.

We first consider concave configurations. The convex and collinear cases will be treated in Section~\ref{sec:global-classification}.

\subsection{Concave kite coordinates}

We first consider the concave case. By \cite[Propositions 1 and 2]{Albouy1995},
every noncollinear planar central configuration of four equal masses possesses
at least one axis of symmetry. Let $q_4$ be the interior body and $q_3$ the
exterior body on the symmetry axis, and let $q_1,q_2$ be the two bodies
symmetric with respect to this axis. Without loss of generality, we assume
\begin{equation}\label{eq:kite-coordinates}
 q_1=(-1,0),\qquad
 q_2=(1,0),\qquad
 q_3=(0,\sqrt{x}),\qquad
 q_4=(0,t\sqrt{x}),
\end{equation}
where
\[
x>0,\qquad 0<t<1.
\]
The mutual distances between the four bodies are
\[
 r_{12}=2,
 \qquad
 r_{13}=r_{23}=\sqrt{1+x},
 \qquad
 r_{14}=r_{24}=\sqrt{1+t^2x},
\]
and 
\[
 r_{34}=\sqrt{x}(1-t).
\]
For the calculations below, we use the notation
\begin{equation} \label{eq:PQR}
 P:=(1+x)^{-\gamma},
 \qquad
 Q:=(1+t^2x)^{-\gamma},
 \qquad
 R:=x^{-\gamma}(1-t)^{1-2\gamma},
 \qquad
 c_A:=2^{1-2\gamma}.
\end{equation}
The center of mass is
\[
 c=\left(0,\frac{(1+t)\sqrt{x}}4\right).
\]
Substituting the above coordinates into \eqref{eq:cc}, the equations at
$q_1$ and $q_2$ give the same two scalar conditions. The horizontal component gives
\begin{equation} \label{eq:lambda}
 \lambda=c_A+P+Q,
\end{equation}
while the vertical component gives
\[
 \frac{\lambda(1+t)}{4}=P+tQ.
\]
At $q_3$ and $q_4$, the horizontal components vanish identically, and
the vertical components give
\begin{equation}\label{eq:vertical-cc}
 \frac{\lambda(3-t)}4=2P+R,
 \qquad
 \frac{\lambda(3t-1)}4=2tQ-R.
\end{equation}
Adding the two equations in \eqref{eq:vertical-cc} yields
\[
 \frac{\lambda(1+t)}4=P+tQ,
\]
which is precisely the vertical equation at $q_1$ and $q_2$. Therefore,
the central-configuration equations \eqref{eq:cc} are equivalent to
\eqref{eq:lambda}--\eqref{eq:vertical-cc}.

Adding and subtracting the two equations in \eqref{eq:vertical-cc}, we obtain
\begin{equation}\label{eq:equivalent}
     \frac{\lambda(1+t)}2=2(P+tQ),
 \qquad
 \lambda(1-t)=2(P-tQ+R).
\end{equation}
Substituting \eqref{eq:lambda} into these two identities gives
\begin{equation} \label{eq:F}
 F(x,t,A):=(3-t)P+(3t-1)Q-c_A(1+t)=0,
\end{equation}
and
\begin{equation} \label{eq:H}
 H(x,t,A):=(1+t)(P-Q)+2R-c_A(1-t)=0.
\end{equation}

Conversely, suppose that $x>0$, $0<t<1$, and \eqref{eq:F}--\eqref{eq:H}
hold. Setting $\lambda=c_A+P+Q$ as in \eqref{eq:lambda},
equations \eqref{eq:F} and \eqref{eq:H} can be rewritten as the two identities in \eqref{eq:equivalent}.
Adding the two identities in
\eqref{eq:equivalent} and dividing by $2$ gives the first equation
in \eqref{eq:vertical-cc}, while subtracting the second identity from
the first and dividing by $2$ gives the second equation in
\eqref{eq:vertical-cc}.
Therefore \eqref{eq:lambda}--\eqref{eq:vertical-cc} hold, and hence
\eqref{eq:cc} is satisfied.

Thus the concave central-configuration problem is equivalent to
\[
 F(x,t,A)=H(x,t,A)=0,
 \qquad x>0,\quad 0<t<1.
\]
It is easy to see that, for every $A\geq 2$,
\begin{equation}\label{eq:triangle-center-point}
 (x,t)=\left(3,\frac13\right)
\end{equation}
is a solution of \eqref{eq:F}--\eqref{eq:H}. It corresponds to the
equilateral triangle with the fourth body at its center. We call \eqref{eq:triangle-center-point} the trivial solution. It remains to determine
whether \eqref{eq:F}--\eqref{eq:H} admit any other solutions.

\section{The degenerate endpoint \texorpdfstring{$A=2$}{A=2}}
\label{sec:endpoint}

In this section, we study the reduced equations at the limiting case
$A=2$. We shall show that the trivial solution
\[
(x,t)=\left(3,\frac13\right)
\]
is the only common zero of $F(\cdot,\cdot,2)$ and
$H(\cdot,\cdot,2)$ in the concave region, and that the corresponding
reduced function has a quadratic zero at $t=1/3$. These properties will
be used in the analysis of the case $2<A\leq3$ in the next section.

The reduced equations \eqref{eq:F}--\eqref{eq:H} are well defined at
$A=2$. Setting $A=2$, we have $\gamma=1$ and
\[
P=\frac1{1+x},\qquad
Q=\frac1{1+t^2x},\qquad
R=\frac1{x(1-t)},\qquad
c_A=\frac12.
\]
Let
\begin{equation}
 w:=tx.
 \label{eq:w-def}
\end{equation}
Then $w>0$ and
\[
 P=\frac{t}{t+w},\qquad
 Q=\frac1{1+tw},\qquad
 R=\frac{t}{w(1-t)}.
\]
Substituting these expressions into \eqref{eq:F}, we obtain
\begin{equation}
 F(x,t,2)
 =-\frac{t(1+t)}{2(t+w)(1+tw)}\,G(t,w),
 \label{eq:F-A2-factor}
\end{equation}
where
\begin{equation}
 G(t,w)
 :=w^2+\left[3\left(t+\frac1t\right)-8\right]w-3.
 \label{eq:G-A2}
\end{equation}
Thus $F(x,t,2)=0$ is equivalent to $G(t,w)=0$. The discriminant of \eqref{eq:G-A2} is positive, and the
product of its two roots is $-3$. Hence \eqref{eq:G-A2} has exactly one
positive root for every $t\in(0,1)$. We denote it by $w_2(t)$ and set
\[
 X_2(t):=\frac{w_2(t)}{t}.
\]
Consequently,
\[
 F(X_2(t),t,2)=0,\qquad 0<t<1.
\]

\begin{proposition}
\label{prop:A2-factorization}
Along the positive root of $F(\cdot,t,2)=0$ one has
\begin{equation}
 H(X_2(t),t,2)
 =\frac{2t^2\bigl(w_2(t)-1\bigr)^2\bigl(w_2(t)+1\bigr)}
 {w_2(t)(1-t)\bigl(t+w_2(t)\bigr)\bigl(1+t w_2(t)\bigr)}.
 \label{eq:h-A2-factor}
\end{equation}
Consequently,
\[
 H(X_2(t),t,2)>0
 \qquad\text{for }t\in(0,1)\setminus\left\{\frac13\right\},
\]
and the only common zero of $F(\cdot,\cdot,2)$ and $H(\cdot,\cdot,2)$ in the concave region is
\[
 (x,t)=\left(3,\frac13\right).
\]
\end{proposition}

\begin{proof}
Substituting \eqref{eq:w-def} into \eqref{eq:H} and simplifying, we obtain
\[
\begin{aligned}
H(x,t,2)
={}&\frac{2t^2(w-1)^2(w+1)}
{w(1-t)(t+w)(1+tw)}\\
&-\frac{t(1+t)^2}
{2(1-t)(t+w)(1+tw)}G(t,w).
\end{aligned}
\]
Since $G(t,w_2(t))=0$, this gives \eqref{eq:h-A2-factor}.

For $0<t<1$, the denominator in \eqref{eq:h-A2-factor}
is positive. Hence $H(X_2(t),t,2)=0$ if and only if
$w_2(t)=1$. Substituting $w=1$ into \eqref{eq:G-A2} yields
\[
3t^2-10t+3=(3t-1)(t-3)=0.
\]
Since $0<t<1$, we obtain $t=1/3$ and $w_2(t)=1$, which
implies $x=w_2(t)/t=3$. This completes the proof.
\end{proof}

Proposition~\ref{prop:A2-factorization} shows that $(3,1/3)$ is the unique common zero
of $F$ and $H$ at $A=2$. We next determine the order of this
zero along the curve $F(x,t,2)=0$.
\begin{proposition}
\label{prop:A2-double-zero}
Let
\[
 h_2(t):=H(X_2(t),t,2).
\]
Then, as $t\to1/3$,
\begin{equation}
 h_2(t)
 =\frac{27}{2}\left(t-\frac13\right)^2
 +O\!\left(\left(t-\frac13\right)^3\right).
 \label{eq:h-A2-expansion}
\end{equation}
In particular,
\[
 h_2\left(\frac13\right)=0,
 \qquad
 h_2'\left(\frac13\right)=0,
 \qquad
 h_2''\left(\frac13\right)=27.
\]
\end{proposition}

\begin{proof}
Differentiating $G(t,w_2(t))=0$ gives
\[
w_2'\left(\frac13\right)
=-\frac{G_t(1/3,1)}{G_w(1/3,1)}=6.
\]
Since $w_2(1/3)=1$, differentiating \eqref{eq:h-A2-factor}
yields
\[
h_2\left(\frac13\right)
=h_2'\left(\frac13\right)=0,
\qquad
h_2''\left(\frac13\right)
=2\cdot\frac38\cdot6^2=27.
\]
The result follows from Taylor's formula.
\end{proof}

\section{Classification of concave configurations for \texorpdfstring{$2<A\le 3$}{2<A<=3}}
\label{sec:concave}

In this section, we study the common zeros of $F$ and $H$ for
$2<A\leq3$. We first show that no common zeros exist for
$0<t<1/3$ or $2/5\leq t<1$. For $1/3\leq t\leq2/5$,
we reduce the problem to a single equation and prove the
strict convexity of the corresponding function. This allows us
to establish the existence and uniqueness of the nontrivial
concave central configuration.

Throughout this section, we assume
\[
 2<A\le3,
 \qquad
 1<\gamma=\frac A2\le\frac32,
\]
and consider the reduced equations \eqref{eq:F}--\eqref{eq:H}.  We shall prove that the permanent solution
\[
 (x,t)=\left(3,\frac13\right)
\]
and one further solution with $1/3<t<2/5$ are the only common zeros in the concave region $x>0$, $0<t<1$.

\subsection{Preliminary facts and the lower interval}

In this subsection, we study the equation $F(x,t,A)=0$
and prove that $F$ and $H$ have no common zeros for
$0<t<1/3$. We first establish the existence and uniqueness
of the positive roots of $F=0$. We then derive estimates
for these roots and use a Jacobian sign argument to exclude
common zeros in the lower interval.

\begin{lemma}
\label{lem:F-root-upper}
For every $t\in[1/3,1)$ the equation $F(x,t,A)=0$ has exactly one positive root $x=X(t,A)$.  Moreover
\[
 F_x(X(t,A),t,A)<0.
\]
\end{lemma}

\begin{proof}
For $t\ge1/3$ and $x>0$, by direct computation, we obtain
\begin{equation}\label{eq:Fx-negative-upper}
    F_x
=-\gamma\left[
\frac{(3-t)P}{1+x}
+\frac{t^2(3t-1)Q}{1+t^2x}
\right]<0.
\end{equation}  Moreover,
\[
 F(0,t,A)=(1+t)(2-c_A)>0,
\]
and
\[
 \lim_{x\to\infty}F(x,t,A)=-c_A(1+t)<0.
\]
Hence $F(x,t,A)=0$ has exactly one positive root.
\end{proof}

 We denote the unique positive root in Lemma~\ref{lem:F-root-upper} by
\[
 x=X(t,A),
\]
and define
\begin{equation}
 h(t,A):=H(X(t,A),t,A).
 \label{eq:h-def}
\end{equation}
We first determine the sign of $h_t$ at $t=1/3$, which will be
used in the subsequent analysis on both sides of this point.

\begin{lemma}
\label{lem:triangle-center-slope}
For every $A>2$,
\begin{equation}
 h_t\left(\frac13,A\right)<0.
 \label{eq:triangle-center-slope}
\end{equation}
\end{lemma}

\begin{proof}
At $t=1/3$ the unique $F$-root is $x=3$.  Differentiating equation \eqref{eq:h-def} directly, we obtain
\begin{equation}
 h_t\left(\frac13,A\right)
 =-\frac92\,4^{-\gamma}\,B(\gamma),
 \qquad
 B(\gamma):=3z^2+(2-4\gamma)z-1,
 \qquad
 z:=3^{\gamma-1}.
 \label{eq:triangle-slope-formula}
\end{equation}
Clearly, $B(1)=0$. Moreover,
\[
 \frac{B'(\gamma)}{2z}
 =(3z-2\gamma+1)\ln3-2.
\]
At $\gamma=1$, we have
\[
 \frac{B'(1)}{2}=2(\ln3-1)>0.
\]
On the other hand,
\[
 \frac{d}{d\gamma}
 \left(\frac{B'(\gamma)}{2z}\right)
 =\ln3\,(3z\ln3-2)>0,
 \qquad \gamma\ge1.
\]
Hence $B'(\gamma)>0$ for $\gamma\ge1$, which implies
$B(\gamma)>0$ for $\gamma>1$. The desired inequality
follows from \eqref{eq:triangle-slope-formula}.
\end{proof}
For $0<t<1/3$, the coefficient $3t-1$ changes sign, so the
monotonicity argument in Lemma~\ref{lem:F-root-upper}
does not apply directly. We next prove that $F=0$ still has
a unique positive root in this interval and establish a
lower bound for this root.

\begin{lemma}
\label{lem:lower-F-root}
For every $t\in[0,1/3)$, the equation $F(u^2,t,A)=0$
has a unique positive root $u=U_-(t,A)$, at which
\[
F_u<0.
\]
Moreover, setting $\alpha:=2^{2/3}-1$, we have
\begin{equation}
U_-(t,A)^2>\alpha>\frac12.
\label{eq:lower-alpha-bound}
\end{equation}
At $t=0$, the equations $F=H=0$ have no common solution.
\end{lemma}

\begin{proof}
Set $m:=A/2+1$ and
\[
D(u):=(3-t)(1+u^2)^{-m}
-t^2(1-3t)(1+t^2u^2)^{-m}.
\]
Differentiating $F(u^2,t,A)$ with respect to $u$, we obtain
\begin{equation}
F_u=-AuD(u).
\label{eq:lower-Fu-factor}
\end{equation}
Let
\[
\phi(u):=
\left(\frac{1+t^2u^2}{1+u^2}\right)^m.
\]
Then
\[
\phi'(u)=
-\frac{2mu(1-t^2)(1+t^2u^2)^{m-1}}
{(1+u^2)^{m+1}}<0,
\qquad u>0.
\]
Since
\[
D(u)=(1+t^2u^2)^{-m}
\bigl[(3-t)\phi(u)-t^2(1-3t)\bigr],
\]
and $D(0)>0$, the sign of $D$ can change at most
once, from positive to negative. Thus $F_u$ can
change sign only from negative to positive.

On the other hand,
\[
F(0,t,A)=(1+t)(2-c_A)>0,
\]
and
\[
\lim_{u\to+\infty}F(u^2,t,A)=
\begin{cases}
-c_A(1+t),&0<t<1/3,\\
-1-c_A,&t=0.
\end{cases}
\]
It follows that $F(u^2,t,A)=0$ has exactly one
positive root, and $F_u<0$ at this root.

Next, let $\alpha=2^{2/3}-1$. Since $3t-1\leq0$ and $Q\leq1$, we have
\[
\begin{aligned}
F(\alpha,t,A)
&\geq (3-t)2^{-A/3}+3t-1-c_A(1+t)\\
&=3\cdot2^{-A/3}-1-c_A
+t(3-2^{-A/3}-c_A)
\geq\frac14.
\end{aligned}
\]
The last inequality follows since
$3\cdot2^{-A/3}-1-c_A$ is decreasing in
$A\in[2,3]$ and equals $1/4$ at $A=3$, while
$3-2^{-A/3}-c_A>0$.
Therefore,
\[
U_-(t,A)^2>\alpha>\frac12.
\]

Finally, if $t=0$, the equations $F=H=0$ imply
$P=u^{-A}$, which contradicts
\[
P=(1+u^2)^{-A/2}<u^{-A}.
\]
\end{proof}

To estimate the Jacobian at a possible common zero of $F$ and $H$,
we first establish an upper bound for the positive root $U_-(t,A)$.

\begin{lemma}
\label{lem:lower-threshold}
For $0\leq t\leq1/3$, let
\begin{equation}
 u_0(t,A)
 :=\sqrt3\left(\frac{2}{3(1-t)}\right)^{1-1/A}.
 \label{eq:lower-u0}
\end{equation}
Then
\begin{equation}
 \frac{d}{dt}F(u_0(t,A)^2,t,A)>0.
 \label{eq:lower-threshold-derivative}
\end{equation}
Consequently, for $0<t<1/3$,
\begin{equation}
 U_-(t,A)<u_0(t,A),
 \qquad U_-(t,A)^2<3.
 \label{eq:lower-u-upper}
\end{equation}
Moreover, setting
\begin{equation}
 z:=2^{A-1}u^{-A}(1-t)^{1-A},
 \qquad
 z_0(A):=3^{A/2-1},
 \label{eq:lower-z}
\end{equation}
we have
\begin{equation}
 z>z_0(A)
 \label{eq:z-threshold}
\end{equation}
at $u=U_-(t,A)$.
\end{lemma}

\begin{proof}
Set
\[
 b:=\frac{2}{3(1-t)},\qquad
 p:=2-\frac1\gamma.
\]
Along \eqref{eq:lower-u0}, we have
\[
 x:=u_0^2=3b^p,
 \qquad
 x'=\frac{(2\gamma-1)x}{\gamma(1-t)}.
\]
Differentiating $F(x(t),t,A)$ with respect to $t$ gives
\begin{equation}
 \frac{dF}{dt}
 =3Q-P-c_A-\mathcal M+\mathcal N,
 \label{eq:lower-threshold-decomp}
\end{equation}
where
\[
 \mathcal M
 =\frac{(3-t)(2\gamma-1)xP}{(1-t)(1+x)},
 \qquad
 \mathcal N=(3t-1)\dot Q\geq0,
\]
and $\dot Q=dQ(x(t),t)/dt\leq0$.

We first estimate $P$. Since $2/3\leq b\leq1$,
\[
 \frac{\partial}{\partial\gamma}
 \ln\left(\gamma b^{2-1/\gamma}\right)
 =\frac{\gamma+\ln b}{\gamma^2}>0.
\]
Hence $\gamma b^p\geq b\geq2/3$. By Bernoulli's
inequality,
\[
 (1+x)^\gamma
 \geq1+\gamma x
 =1+3\gamma b^p\geq3,
\]
which gives $P\leq1/3$.

Next, let
\[
 f(b):=\gamma(1+3b^p)-3b(2\gamma-1).
\]
Since
\[
 f'(b)=3(2\gamma-1)(b^{p-1}-1)\leq0,
 \qquad
 f(1)=3-2\gamma\geq0,
\]
we obtain $3b(2\gamma-1)\leq\gamma(1+x)$.
Using $(3-t)/(1-t)=3b+1$, we have
\[
 \mathcal M
 \leq\frac{x(\gamma x+3\gamma-1)}
 {(1+x)^{\gamma+1}}.
\]
Using Taylor's formula for $(1+x)^{\gamma+1}$, we have
\[
\begin{aligned}
 &(1+x)^{\gamma+1}
 -x(\gamma x+3\gamma-1)\\
 &\quad\geq
 1-2(\gamma-1)x
 +\frac{\gamma(\gamma-1)}2x^2\\
 &\quad=
 \frac{\gamma(\gamma-1)}2
 \left(x-\frac2\gamma\right)^2
 +\frac{2-\gamma}{\gamma}>0.
\end{aligned}
\]
Therefore $\mathcal M<1$.

Since $x\leq3$ and $t\leq1/3$, we also have
\[
 Q\geq\left(\frac43\right)^{-3/2},
 \qquad c_A\leq\frac12.
\]
Combining these estimates with
\eqref{eq:lower-threshold-decomp}, we obtain
\[
 \frac{dF}{dt}
 >3\left(\frac43\right)^{-3/2}
 -\frac13-\frac12-1>0.
\]
This proves \eqref{eq:lower-threshold-derivative}.

At $t=1/3$, we have $u_0=\sqrt3$ and
$F(3,1/3,A)=0$. Thus
\[
 F(u_0(t,A)^2,t,A)<0,
 \qquad 0\leq t<\frac13.
\]
Since $F(0,t,A)>0$ and $F(u^2,t,A)=0$ has a
unique positive root, it follows that
\[
 U_-(t,A)<u_0(t,A).
\]
Moreover, $u_0(t,A)^2=3b^p<3$ for $t<1/3$,
which proves \eqref{eq:lower-u-upper}.

Finally, $z$ is strictly decreasing in $u$, and
direct substitution gives
\[
 z(u_0(t,A),t,A)=z_0(A).
\]
Hence $U_-(t,A)<u_0(t,A)$ implies
$z>z_0(A)$.
\end{proof}

By Lemmas~\ref{lem:lower-F-root} and
\ref{lem:lower-threshold}, any common zero of $F$ and $H$
with $0<t<1/3$ must satisfy
$\frac12<u^2<3$ and $z>z_0(A)$. We now use these estimates to determine the sign of the
Jacobian $F_uH_t-F_tH_u$ at such a point.

\begin{lemma}
\label{lem:lower-jacobian}
Let $0<t<1/3$, and suppose $F=H=0$.  Then
\begin{equation}
 \Delta:=F_uH_t-F_tH_u>0.
 \label{eq:lower-Delta-positive}
\end{equation}
\end{lemma}

\begin{proof}
Write $x=u^2$ and let $z$ be defined by \eqref{eq:lower-z}. Since $R=c_Az$, solving $F=H=0$ for $P$ and $Q$ gives
\begin{equation}
 P=c_A\frac{t(3-t)-(3t-1)z}{(1+t)^2},
 \qquad
 Q=c_A\frac{(3t-1)+(3-t)z}{(1+t)^2}.
 \label{eq:lower-PQ-elimination}
\end{equation}
Substituting these expressions into the Jacobian, we obtain
\begin{equation}
 \Delta
 =\frac{2Ac_A^2}{u(1-t)(1+t)^2(1+u^2)(1+t^2u^2)}\,
 \mathcal P(A,t,x,z),
 \label{eq:lower-Delta-factor}
\end{equation}
where
\begin{align}
 \mathcal R
 &:=(1-t)(3-t)(3t-1)(xt^2-z^2)
    +2t(5t^2-6t+5)z(xt+1),
 \label{eq:lower-Rcal}\\
 \mathcal S
 &:=x^2t^2(3t^2-2t+3)
    +9xz(1-t)^3+x t(15t^2-26t+15)+4z(1-t),
 \label{eq:lower-Scal}\\
 \mathcal P
 &:=-Ax\mathcal R+(t+2z-1)\mathcal S.
 \label{eq:lower-Pcal}
\end{align}
The coefficient of $\mathcal P$ in
\eqref{eq:lower-Delta-factor} is positive.
By Lemmas~\ref{lem:lower-F-root} and
\ref{lem:lower-threshold}, it is enough to prove
$\mathcal P>0$ for
\begin{equation}
2\leq A\leq3,\quad
0\leq t\leq\frac13,\quad
\frac12\leq x\leq3,\quad
z>z_0(A)\geq1.
\label{eq:lower-jacobian-domain}
\end{equation}
Set
\[
\varphi(t):=5t^2-6t+5,\qquad
\chi(t):=3t^2-10t+3.
\]
Both functions are nonnegative for
$0\leq t\leq1/3$, and $\varphi(t)>0$.
Direct differentiation gives
\[
\begin{aligned}
\frac12\mathcal P_{zz}
&=(1-t)\bigl[
x\bigl(A(3-t)(3t-1)+18(1-t)^2\bigr)+8
\bigr]>0,\\
-\frac1{2x}\partial_A\mathcal P_z
&=t^2x\varphi(t)
+z(1-t)\chi(t)+t\varphi(t)>0.
\end{aligned}
\]
Thus $\mathcal P_z$ is increasing in $z$ and
decreasing in $A$. It is therefore sufficient
to consider
\[
G(t,x):=\mathcal P_z(3,t,x,1).
\]
Direct computations give
\[
\begin{aligned}
G_{xx}
&=-16t^2(3t^2-4t+3)\leq0,\\
G\left(t,\frac12\right)
&\geq\frac{33-10t-60t^2}{2}>0,\\
G(t,3)
&\geq39+10t-340t^2>0.
\end{aligned}
\]
By the concavity of $G$ in $x$, we obtain $G(t,x)>0$
for $0\leq t\leq1/3$ and $1/2\leq x\leq3$.
we conclude that
\begin{equation}
\mathcal P_z>0
\label{eq:lower-Pz-positive}
\end{equation}
throughout the domain
\eqref{eq:lower-jacobian-domain}.

Now define
\begin{equation}
D_0(A,t,x,z):=
-x\mathcal R+\frac z2\mathcal P_z.
\label{eq:lower-D0}
\end{equation}
We have
\[
-\frac1x\partial_A(D_0)_z
=t^2x\varphi(t)+2z(1-t)\chi(t)
+t\varphi(t)>0.
\]
Moreover, $(D_0)_{zz}$ is decreasing in $A$,
and
\[
(D_0)_{zz}\big|_{A=3}
=4(1-t)(3t^2x+2tx+3x+4)>0.
\]
Hence $(D_0)_z$ is decreasing in $A$ and
increasing in $z$.

At $A=3$ and $z=1$, direct computation gives
\[
\partial_{xx}(D_0)_z
=-4t^2(11t^2-14t+11)\leq0.
\]
Furthermore,
\[
\begin{aligned}
(D_0)_z\left(3,t,\frac12,1\right)
&\geq\frac{71-40t-68t^2-31t^4}{4}>0,\\
(D_0)_z(3,t,3,1)
&\geq\frac{73-514t^2}{2}>0.
\end{aligned}
\]
By concavity in $x$, we obtain
\begin{equation}
(D_0)_z>0.
\label{eq:lower-D0z-positive}
\end{equation}

Since $\ln3>1$, we have $z_0(A)\geq \frac A2$. Set $E(A,t,x):=D_0\left(A,t,x,\frac A2\right)$.
We claim that $E(A,t,x)\geq0$ on the
domain under consideration. To prove this, we first show that $E_A$ is
positive. Direct differentiation gives
\[
(E_A)_{AA}
=-\frac32x(3-t)(1-t)(1-3t)\leq0.
\]
In addition,
\[
-2(E_A)_{tt}=a(A,t)x^2+b(A,t)x+4>0,
\]
where
\[
\begin{aligned}
a(A,t)
&=A(120t^2-72t+20)
  +84t^2-60t+14,\\
b(A,t)
&=A^2(39-27t)+A(150t-106)\\
&\quad+54t^2-138t+82.
\end{aligned}
\]
Both coefficients are positive for
$2\leq A\leq3$ and $0\leq t\leq1/3$.
For $b$, this follows from
\[
b(2,t)=26+54t+54t^2>0,
\qquad
\partial_A b>0.
\]

Therefore $E_A$ is concave in both $A$ and
$t$. Its minimum is bounded below by its
values at the four corners of
$[2,3]\times[0,1/3]$. These values are
\[
\begin{aligned}
E_A(2,0,x)&=\frac{15x+28}{4}>0,\quad
&E_A(3,0,x)&=11>0,\\
E_A(2,\tfrac13,x)&=\tfrac4{27}(-7x^2+18x+33)>0,\quad
&E_A(3,\tfrac13,x)&=\tfrac4{81}(-29x^2+84x+153)>0.
\end{aligned}
\]
for $1/2\leq x\leq3$. Consequently,
\begin{equation}
E_A>0.
\label{eq:lower-EA-positive}
\end{equation}

At $A=2$, the function $E(2,t,x)$ is
concave in $x$, and
\[
\begin{aligned}
E\left(2,t,\frac12\right)
&=\frac{
33-20t-50t^2+17t^3-13t^4-3t^5}{4}>0,\\
E(2,t,3)
&=\frac{1-9t^2}{2}
(6t^3+11t^2-20t+39)\geq0.
\end{aligned}
\]
It follows from \eqref{eq:lower-EA-positive}
that
\[
E(A,t,x)\geq E(2,t,x)\geq0.
\]
Together with $z_0(A)\geq A/2$ and
\eqref{eq:lower-D0z-positive}, this yields
\[
D_0(A,t,x,z_0(A))\geq0.
\]

We now differentiate $\mathcal P$ along
$z=z_0(A)$. By \eqref{eq:lower-D0},
\[
\begin{aligned}
\frac{d}{dA}\mathcal P(A,t,x,z_0(A))
&=-x\mathcal R+z_0'(A)\mathcal P_z\\
&=D_0+
\frac{\ln3-1}{2}z_0(A)\mathcal P_z>0.
\end{aligned}
\]
Thus $\mathcal P(A,t,x,z_0(A))$ is
increasing in $A$.

At $A=2$, we have $z_0(2)=1$. Set
\[
Q_0(t,x):=\frac{\mathcal P(2,t,x,1)}{1+t}.
\]
The function $Q_0$ is concave in $x$,
and
\[
\begin{aligned}
Q_0\left(t,\frac12\right)
&=\frac{22-16t-37t^2+22t^3-3t^4}{4}>0,\\
Q_0(t,3)
&=(1-3t)(9t^3-33t^2+35t+13)\geq0.
\end{aligned}
\]
Hence
\[
\mathcal P(A,t,x,z_0(A))\geq\mathcal P(2,t,x,1)\geq0.
\]

Finally, the common zero satisfies
$z>z_0(A)$. Since $\mathcal P_z>0$,
we conclude that
\[
\mathcal P(A,t,x,z)>
\mathcal P(A,t,x,z_0(A))\geq0.
\]
The result follows from
\eqref{eq:lower-Delta-factor}.
\end{proof}
We now use the Jacobian estimate in
Lemma~\ref{lem:lower-jacobian} to show that
$F$ and $H$ have no common zeros for $0<t<1/3$.
\begin{proposition}
\label{prop:no-lower-root}
There is no common zero of $F$ and $H$ with $0\le t<1/3$.
\end{proposition}

\begin{proof}
The case $t=0$ has been excluded in
Lemma~\ref{lem:lower-F-root}. For $0<t<1/3$, let
$u=U_-(t,A)$ be the unique positive root of $F=0$
and define
\[
h_-(t,A):=H(U_-(t,A)^2,t,A).
\]
Since $F_u<0$, the implicit function theorem
implies that $U_-$ is continuously differentiable. By Lemma~\ref{lem:lower-jacobian}, at any zero of $h_-$, we have
\[
(h_-)_t=\frac{\Delta}{F_u}<0.
\]
Moreover, $h_-$ extends smoothly to $t=1/3$, where
\[
h_-\left(\frac13,A\right)=0,
\qquad
(h_-)_t\left(\frac13,A\right)<0.
\]
Hence $h_-(t,A)>0$ for $t<1/3$ sufficiently
close to $1/3$.

Suppose that $h_-(t_0,A)=0$ for some
$t_0\in(0,1/3)$. Since
$(h_-)_t(t_0,A)<0$, we have
$h_-(t,A)<0$ for $t>t_0$ sufficiently
close to $t_0$.
On the other hand, $h_-(t,A)>0$ for
$t<1/3$ sufficiently close to $1/3$.
By continuity, there must be another zero
where $h_-$ changes from negative to
positive. At such a zero, its derivative
cannot be negative, contradicting
$(h_-)_t<0$ at every zero.
\end{proof}

\subsection{Exclusion of the upper interval \texorpdfstring{$t\ge2/5$}{t>=2/5}}

We now show that $H$ is positive at the unique positive
root of $F=0$ when $t\geq2/5$.

\begin{proposition}
\label{prop:upper-exclusion}
For every $2\leq A\leq3$ and $2/5\leq t<1$,
if $F(x,t,A)=0$ with $x>0$, then
\[
 H(x,t,A)>0.
\]
\end{proposition}

\begin{proof}
Using $F=0$ to eliminate $c_A$ from \eqref{eq:H},
we obtain
\begin{equation}
 H=\frac{2Q}{1+t}L_\gamma,
 \label{eq:H-upper-L}
\end{equation}
where $L_\gamma$ is defined by
\begin{equation}
 L_\gamma
 =(3t-1)a^\gamma+(1-t^2)b^\gamma-t(3-t),
 \label{eq:Lgamma}
\end{equation}
and $a$ and $b$ are given by
\[
 a:=\frac{1+t^2x}{1+x},
 \qquad
 b:=\frac{1+t^2x}{x(1-t)^2}.
\]
For fixed $t$ and $x$, the function $L_\gamma$ is convex
in $\gamma$. Applying the inequality $y\ln y\geq y-1$
for $y>0$, we obtain
\[
\begin{aligned}
 \left.\partial_\gamma L_\gamma\right|_{\gamma=1}
 &=(3t-1)a\ln a+(1-t^2)b\ln b\\
 &\geq(3t-1)(a-1)+(1-t^2)(b-1)
 =L_1.
\end{aligned}
\]
 Hence, by convexity, $L_\gamma\geq\gamma L_1$
for $\gamma\geq1$, which implies
\begin{equation}
 L_1>0\quad\Longrightarrow\quad L_\gamma>0.
 \label{eq:exponent-lifting}
\end{equation}

For $\gamma=1$, set $\Pi(t,x):=t(3t^2-7t+3)x^2-2tx-1$. Direct computation gives
\begin{equation}
 L_1
 =\frac{1+t}{x(t-1)(x+1)}\Pi(t,x).
 \label{eq:L1-Pi}
\end{equation}
Since the coefficient in \eqref{eq:L1-Pi} is
negative, it remains to prove that
$\Pi(t,X(t,A))<0$.

Set
\[
 \tau:=\frac{7-\sqrt{13}}6.
\]
For $\tau\leq t<1$, we have
$3t^2-7t+3\leq0$, and hence
$\Pi(t,x)<0$ for every $x>0$.

Now consider $2/5\leq t<\tau$.
In this interval, $\Pi(t,\cdot)$ has
a unique positive root $r(t)$. Writing
$s=\sqrt{3t}$, we obtain
\[
 r(t)=-\frac3{s(s^2+s-3)},
 \qquad
 \frac{dr}{ds}
 =\frac{3(3s^2+2s-3)}
 {s^2(s^2+s-3)^2}>0.
\]
Thus $r(t)$ is strictly increasing.

We next compare $X(t,A)$ with $r(t)$.
Define
\begin{equation}
\begin{aligned}
 \Phi(\gamma,t,x)
 :=\frac{F(x,t,A)}{c_A}
 ={}&\frac{3-t}{2}
 \left(\frac4{1+x}\right)^\gamma
 +\frac{3t-1}{2}
 \left(\frac4{1+t^2x}\right)^\gamma
 -(1+t).
\end{aligned}
\label{eq:Phi-upper}
\end{equation}

For $2/5\leq t\leq9/20$ and
$3\leq x\leq9/2$, differentiation gives
\[
\begin{aligned}
 \Phi_t
 ={}&-\frac12\left(\frac4{1+x}\right)^\gamma-1
 +\left(\frac32-
 \frac{\gamma tx(3t-1)}{1+t^2x}\right)
 \left(\frac4{1+t^2x}\right)^\gamma.
\end{aligned}
\]
On this domain,
\[
 \left(\frac4{1+x}\right)^\gamma\leq1,
 \qquad
 \left(\frac4{1+t^2x}\right)^\gamma>2,
 \qquad
 \frac{\gamma tx(3t-1)}{1+t^2x}<\frac34.
\]
Consequently,
\begin{equation}
 \Phi_t>-\frac12-1+\frac34\cdot2=0.
 \label{eq:Phi-t-lower}
\end{equation}

For fixed $(t,x)$, the function
$\gamma\mapsto\Phi(\gamma,t,x)$ is convex.
The following table gives four intervals
$[t_-,t_+]$ and corresponding values $x_*$:
\[
\begin{array}{c|c|c}
t_- & t_+ & x_*\\ \hline
{2}/{5} & {41}/{100} & {31}/{8}\\
{41}/{100} & {21}/{50} & 4\\
{21}/{50} & {11}/{25} & {17}/{4}\\
{11}/{25} & {9}/{20} & {9}/{2}
\end{array}
\]
For each row, direct substitution gives
\[
 \Pi(t_-,x_*)<0,
 \qquad
 \Phi(1,t_+,x_*)<0,
 \qquad
 \Phi(3/2,t_+,x_*)<0.
\]
These inequalities can be verified exactly;
the square roots at $\gamma=3/2$ are
handled by rational comparisons.

Since $\Phi_t>0$ and $\Phi$ is convex
in $\gamma$, we obtain
\[
 \Phi(\gamma,t,x_*)<0,
 \qquad
 t_-\leq t\leq t_+,\quad
 1\leq\gamma\leq\frac32.
\]
As $F_x<0$, this implies
\[
 X(t,A)<x_*.
\]
Moreover, $\Pi(t_-,x_*)<0$ and the
monotonicity of $r$ give
\[
 x_*<r(t_-)\leq r(t).
\]
Hence
\[
 \Pi(t,X(t,A))<0,
 \qquad \frac25\leq t\leq\frac9{20}.
\]

For $9/20\leq t<\tau$, we use
$x_*=5$. At $\gamma=1$,
\[
 \Phi(1,t,5)
 =-\frac{2(t+1)(10t^2-10t+3)}
 {3(5t^2+1)}<0.
\]
At $\gamma=3/2$, we have
\[
 \Phi(3/2,t,5)
 =\frac{3-t}{2}\left(\frac23\right)^{3/2}
 +\frac{4(3t-1)}{(1+5t^2)^{3/2}}
 -(1+t).
\]
Since $\tau<3/5$ and
\[
 \frac{d}{dt}
 \frac{4(3t-1)}{(1+5t^2)^{3/2}}
 =\frac{12(1+5t-10t^2)}
 {(1+5t^2)^{5/2}}>0
\]
for $9/20\leq t\leq3/5$, we obtain
\[
\begin{aligned}
 \frac{3-t}{2}\left(\frac23\right)^{3/2}
 \leq\frac{17}{10\sqrt6}<\frac7{10}\quad \text{and} \quad
 \frac{4(3t-1)}{(1+5t^2)^{3/2}}
 \leq\sqrt{\frac{160}{343}}<\frac7{10}.
\end{aligned}
\]
Therefore
\[
 \Phi(3/2,t,5)
 <\frac75-(1+t)<0.
\]
By convexity in $\gamma$, we conclude that
$\Phi(\gamma,t,5)<0$ for $1\leq\gamma\leq3/2$.
Since $F_x<0$, we get $X(t,A)<5$.

Since $\Pi(9/20,5)<0$ and $r$ is increasing, we obtain
$X(t,A)<5<r(9/20)\leq r(t)$ for $9/20\leq t<\tau$.
Hence $\Pi(t,X(t,A))<0$ on this interval.
We have proved $\Pi(t,X(t,A))<0$
for every $2/5\leq t<1$. Therefore
$L_1>0$, and \eqref{eq:exponent-lifting}
gives $L_\gamma>0$. By \eqref{eq:H-upper-L},
we conclude that $H(x,t,A)>0$.
\end{proof}

\subsection{Strict convexity on the remaining interval}

We now consider the remaining interval
$1/3\leq t\leq2/5$. By Lemma~\ref{lem:F-root-upper},
the equation $F(x,t,A)=0$ has a unique positive
root $x=X(t,A)$. We first establish the bounds
needed in the subsequent analysis.

\begin{lemma}
\label{lem:upper-branch-box}
For $1/3\le t\le2/5$ and $2\le A\le3$, the positive root $X(t,A)$ satisfies
\begin{equation}
 3\leq X(t,A)<4.
 \label{eq:X-box}
\end{equation}
Moreover, let
\begin{equation}
 Y:=\frac{P}{Q}
 =\left(\frac{1+t^2x}{1+x}\right)^\gamma,
 \qquad
 Z:=\frac{R}{Q}
 =(t^2+x^{-1})^\gamma(1-t)^{1-2\gamma},
 \label{eq:YZ-def}
\end{equation}
then, at $x=X(t,A)$,
\begin{equation}
 \frac3{20}\le Y\le\frac38,
 \qquad
 \frac{61}{125}<Z<1.
 \label{eq:YZ-box}
\end{equation}
\end{lemma}

\begin{proof}
At $x=3$, we have
\[
 F(3,t,A)
 =(3t-1)\left[(1+3t^2)^{-\gamma}-4^{-\gamma}\right]\ge0.
\]
Since $F_x<0$, it follows that $X(t,A)\geq3$. 

 To prove $F(4,t,A)<0$, set
\[
 G(\gamma,t):=4^\gamma F(4,t,A)
 =(3-t)\left(\frac45\right)^\gamma
 +(3t-1)\left(\frac4{1+4t^2}\right)^\gamma-2(1+t).
\]
For fixed $t$, the function $G$ is convex
in $\gamma$. It is therefore sufficient to
prove $G(1,t)<0$ and $G(3/2,t)<0$. At $\gamma=1$, direct computation gives
\[
 G(1,t)
 =-\frac{2(t+1)(2t-1)(14t-9)}{5(4t^2+1)}<0.
\]
At $\gamma=3/2$, we have
\[
 G_t
 =24\frac{1+4t-8t^2}{(1+4t^2)^{5/2}}-2-\frac8{5\sqrt5}>0
\]
on $[1/3,2/5]$, while
\[
 G\left(\frac32,\frac25\right)
 =\frac{104}{25\sqrt5}+\frac{200}{41\sqrt{41}}-\frac{14}{5}
<0.
\]
Thus $G(3/2,t)<0$, and convexity in
$\gamma$ gives $G(\gamma,t)<0$ for
$1\leq\gamma\leq3/2$. Hence
$F(4,t,A)<0$, which implies
$X(t,A)<4$.

Next, the expression $(1+t^2x)/(1+x)$ is increasing in $t$ and decreasing
in $x$. Therefore,
\[
 \frac{13}{45}\le a\le\frac{37}{100}.
\]
Since $1\leq\gamma\leq3/2$, we obtain
\[
 \left(\frac{13}{45}\right)^{3/2}
 \leq Y\leq\frac{37}{100},
\]
and hence $\frac3{20}\leq Y\leq\frac38$.

For $Z$, taking logarithmic differentiation yields
\[
 (\ln Z)_t
 =\frac{2\gamma t}{t^2+x^{-1}}+\frac{2\gamma-1}{1-t}>0,
 \qquad
 (\ln Z)_x
 =-\frac{\gamma}{x^2(t^2+x^{-1})}<0.
\]
Thus, for fixed $\gamma$, the minimum and maximum occur at $(t,x)=(1/3,4)$ and $(2/5,3)$, respectively.  To determine the dependence of $Z$ on $\gamma$, we write
\[
Z=(1-t)\left[\frac{t^2+x^{-1}}{(1-t)^2}\right]^\gamma.
\]
At $(t,x)=(1/3,4)$ and $(2/5,3)$, the expressions
in brackets are $13/16<1$ and $37/27>1$, respectively.
Hence both extremal values over $1\leq\gamma\leq3/2$
occur at $\gamma=3/2$. It follows that
\[
 Z\ge\frac{13\sqrt{13}}{96}>\frac{61}{125},
 \qquad
 Z\le\frac35\left(\frac{37}{27}\right)^{3/2}<1.
\]
This proves \eqref{eq:YZ-box}.
\end{proof}

We next study the convexity of $h(t,A)$ on $1/3\leq t\leq2/5$. We first derive an identity for $h_{tt}$ using the relation $F(X(t,A),t,A)=0$.
Throughout the following calculations, $A$ is fixed, and all derivatives are taken with respect to $x$ and $t$. The relation $F=0$ will be used only
after differentiation.

\begin{lemma}
\label{lem:relative-Hessian}
Let $x=X(t)$ be a smooth function satisfying $F(X(t),t)=0$ and $F_x\ne0$. Set $h(t)=H(X(t),t)$ and define
\begin{equation}
 \mathcal K:=H_xD^2F-F_xD^2H
 =\begin{pmatrix}K_{11}&K_{12}\\K_{12}&K_{22}\end{pmatrix}.
 \label{eq:relative-Hessian}
\end{equation}
Then
\begin{equation}
 (-F_x)h_{tt}
 =(X_t,1)\mathcal K\binom{X_t}{1}
 =K_{11}X_t^2+2K_{12}X_t+K_{22}.
 \label{eq:relative-Hessian-identity}
\end{equation}
\end{lemma}

\begin{proof}
Differentiating $F(X(t),t)=0$ with respect to $t$, we obtain
\[
 X_t=-\frac{F_t}{F_x},
 \qquad
 X_{tt}
 =-\frac{F_{xx}X_t^2+2F_{xt}X_t+F_{tt}}{F_x}.
\]
On the other hand, 
\[
 h_{tt}=H_{xx}X_t^2+2H_{xt}X_t+H_{tt}+H_xX_{tt}
\]
Substituting the expression for $X_{tt}$ gives \eqref{eq:relative-Hessian-identity}.
\end{proof}

We now apply Lemma~\ref{lem:relative-Hessian} to the functions $F$ and $H$. After differentiation and using $F=0$ to eliminate $c_A$, we obtain
\begin{equation}
\begin{aligned}
 K_{11}=2\gamma^2(\gamma+1)\Bigg[
 \frac{t^2(1-t)(1+t)^3PQ}{(1+x)^2(1+t^2x)^2}
 +\frac{(3-t)PR}{x^2(1+x)^2}
 +\frac{t^2(3t-1)QR}{x^2(1+t^2x)^2}
 \Bigg].
\end{aligned}
\label{eq:K11-positive}
\end{equation}
Since $P,Q,R>0$ and $1/3\leq t\leq2/5$, we have $K_{11}>0$. Together with $F_x<0$ and \eqref{eq:relative-Hessian-identity}, it suffices to prove
\[
K_{22}+2K_{12}X_t>0.
\]
This is established in the following proposition.

\begin{lemma}
\label{prop:direct-tangent-positive}
For $1/3\leq t\leq2/5$ and $1\leq\gamma\leq3/2$, the unique positive root $x=X(t,A)$ of $F=0$ satisfies
\begin{equation}
 K_{22}+2K_{12}X_t>0.
 \label{eq:direct-tangent-positive}
\end{equation}
\end{lemma}

\begin{proof}
Since $X_t=-F_t/F_x$, we have 
\begin{equation}
 K_{22}+2K_{12}X_t
 =\frac{K_{22}F_x-2K_{12}F_t}{F_x}.
 \label{eq:direct-tangent-ratio}
\end{equation}
Thus, since $F_x<0$, it is sufficient to prove $K_{22}F_x-2K_{12}F_t<0$.
We first differentiate $F$ and $H$ with $c_A$ fixed. Restricting to $F=0$, we substitute
\[
 c_A=\frac{(3-t)P+(3t-1)Q}{1+t}.
\]
Using $Y=P/Q$ and $Z=R/Q$, and clearing
denominators, we obtain
\begin{equation}
 -\frac{K_{22}F_x-2K_{12}F_t}{4\gamma^2Q^3}
 =\frac{\mathscr M(\gamma,t,x,Y,Z)}
 {x(1-t)^2(1+t)(1+x)^2(1+t^2x)^3},
 \label{eq:M-identity}
\end{equation}
Here $\mathscr M$ is a polynomial in $(\gamma,t,x,Y,Z)$ with rational coefficients. Its exact reconstruction is given in Appendix~\ref{app:M-reconstruction}. Since the denominator in \eqref{eq:M-identity}
is positive, it remains to prove
$\mathscr M>0$. By Lemma~\ref{lem:upper-branch-box},
we have
\[
\begin{gathered}
1\leq\gamma\leq\frac32,\qquad
\frac13\leq t\leq\frac25,\qquad
3\leq x\leq4,\\
Y\geq Y_0:=\frac3{20},
\qquad
Z>Z_0:=\frac{61}{125}.
\end{gathered}
\]

We first consider the dependence of
$\mathscr M$ on $Z$. Differentiation gives
\[
 \mathscr M_Z=aY^2+bY+d,
\]
where $a,b,d$ depend on $(\gamma,t,x)$.
The estimates in Appendix~\ref{app:Z-certificate} give
\[
 a>0,\qquad 4ad-b^2>0.
\]
The bounds for $a$ and $d$ follow from
factorization and monotonicity, while
the estimate for $b$ uses the fifteen
exact coefficient inequalities in Table~\ref{tab:app-bt-certificate}
of Appendix~\ref{app:exact-algebra}. Consequently, $\mathscr M_Z>0$ for every real $Y$.

It remains to estimate $\mathscr M$ at $Z=Z_0$. Write
\[
 \mathscr M=A_3Y^3+B_2Y^2+C_1Y+C_0
 \quad \text{where} \quad
 A_3=8x(1-t)^2(1+t^2x)^3>0.
\]
The exact polynomial estimates in
Appendix~\ref{app:Y-certificate} establish
\begin{equation}
\begin{gathered}
 B_2>0,\quad
 J:=2Y_0B_2+C_1>0,\quad\text{and}\quad
 \mathscr M(\gamma,t,x,Y_0,Z_0)>0.
\end{gathered}
\label{eq:M-coefficient-positive}
\end{equation}
These inequalities are verified by the shifted-power coefficient method described in Appendix~\ref{app:coefficient-principle}, with the exact certificates listed in Table~\ref{tab:app-shifted-certificates}.

For $Y\geq Y_0$, it follows that
\[
\begin{aligned}
 \mathscr M_Y(\gamma,t,x,Y,Z_0)
 =3A_3Y^2+2B_2Y+C_1
 =3A_3Y^2+2B_2(Y-Y_0)+J>0.
\end{aligned}
\]
Hence $\mathscr M$ is increasing in $Y$
when $Z=Z_0$. Combining this with
$\mathscr M_Z>0$ and
\eqref{eq:M-coefficient-positive}, we obtain
\[
\begin{aligned}
 \mathscr M(\gamma,t,x,Y,Z)
 >\mathscr M(\gamma,t,x,Y,Z_0)
 \geq\mathscr M(\gamma,t,x,Y_0,Z_0)>0.
\end{aligned}
\]

Therefore \eqref{eq:M-identity} implies
$K_{22}F_x-2K_{12}F_t<0$. Since $F_x<0$,
equation \eqref{eq:direct-tangent-ratio}
gives \eqref{eq:direct-tangent-positive}.
\end{proof}
Combining the preceding estimates with Lemma~\ref{lem:relative-Hessian}, we obtain the strict convexity of the reduced function.
\begin{theorem}
\label{thm:strict-convexity}
For every $2\le A\le3$, the reduced function $h(t,A)$ satisfies
\begin{equation}
 h_{tt}(t,A)>0,
 \qquad
 \frac13\le t\le\frac25.
 \label{eq:h-strict-convexity}
\end{equation}
\end{theorem}

\begin{proof}
By Lemma~\ref{lem:relative-Hessian}, Lemma~\ref{prop:direct-tangent-positive}, and \eqref{eq:K11-positive}, we have
\[
 (-F_x)h_{tt}
 =K_{11}X_t^2+\bigl(K_{22}+2K_{12}X_t\bigr)>0.
\]
Since $F_x<0$, the result follows.
\end{proof}

\subsection{The unique nontrivial concave branch}

By Propositions~\ref{prop:no-lower-root} and~\ref{prop:upper-exclusion}, any nontrivial common zero of $F$ and $H$ must lie in $1/3<t<2/5$. Together with the strict convexity of $h(t,A)$ established in Theorem~\ref{thm:strict-convexity}, we obtain the following classification result.

\begin{theorem}
\label{thm:concave-classification}
For each $2<A\leq3$, the system $F=H=0$
admits two solutions within $x>0,\,0<t<1$: the trivial solution $(3,1/3)$ and
a unique nontrivial solution 
\begin{equation}
 \bigl(X(t_*(A),A),t_*(A)\bigr),
 \qquad
 \frac13<t_*(A)<\frac25.
 \label{eq:nontrivial-concave-root}
\end{equation}
At this nontrivial solution,
\begin{equation}
 h_t(t_*(A),A)>0.
 \label{eq:positive-upper-crossing}
\end{equation}
\end{theorem}

\begin{proof}
By Propositions~\ref{prop:no-lower-root} and \ref{prop:upper-exclusion}, there are no common zeros of $F$ and $H$ for $0<t<1/3$ or $t\geq2/5$.
Thus it remains to consider $1/3\leq t\leq2/5$. 

By Theorem~\ref{thm:strict-convexity}, the function
$h(\cdot,A)$ is strictly convex on $[1/3,2/5]$.
Moreover, Lemma~\ref{lem:triangle-center-slope}
and Proposition~\ref{prop:upper-exclusion} give
\[
h\left(\frac13,A\right)=0,\qquad
h_t\left(\frac13,A\right)<0,\qquad
h\left(\frac25,A\right)>0.
\]
Hence $h(t,A)<0$ for $t>1/3$ sufficiently close to $1/3$. By continuity, there exists $t_*(A)\in(1/3,2/5)$ such that $h(t_*(A),A)=0$. To prove uniqueness, suppose that $h(t_1,A)=h(t_2,A)=0$ for some $1/3<t_1<t_2<2/5$. Since $h$ is strictly convex and $h(1/3,A)=h(t_2,A)=0$, we obtain $h(t_1,A)<0$, which is a contradiction.

Finally, if $h_t(t_*,A)\le0$, the strict monotonicity of $h_t$ implies $h_t<0$ on the entire interval $(1/3,t_*)$. This contradicts
\[
 h(t_*,A)-h\left(\frac13,A\right)=0.
\]
Thus \eqref{eq:positive-upper-crossing} holds.
\end{proof}

Theorem~\ref{thm:concave-classification} settles the only geometric class in which a nontrivial one-axis configuration can occur.  In the next section we exclude all additional convex and degenerate noncollinear types and thereby complete the planar classification for $2<A\le3$.

\section{Completion of the planar classification}
\label{sec:global-classification}

In this section, we consider the remaining cases and complete the classification for $2<A\leq3$. We first exclude noncollinear configurations with exactly three collinear bodies. We then study the convex configurations and finally consider the collinear case.

\begin{lemma}
\label{lem:no-three-collinear}
A noncollinear central configuration of four equal masses cannot have exactly three collinear bodies.
\end{lemma}

\begin{proof}
Suppose that $q_1,q_2,q_3$ are collinear, while $q_4$ does not lie on the same line. Without loss of generality, we assume
\[
q_i=(x_i,0),\quad i=1,2,3,
\qquad q_4=(a,b),\quad b\ne0,
\]
where $x_1,x_2,x_3$ are distinct. Since the vertical component of the center of mass is $b/4$, the vertical component of \eqref{eq:cc} gives
\[
-\frac{\lambda b}{4}=-\frac{b}{r_{i4}^A},
\qquad i=1,2,3.
\]
Hence $r_{14}=r_{24}=r_{34}$, which implies
\[
(x_1-a)^2=(x_2-a)^2=(x_3-a)^2.
\]
This is impossible since $x_1,x_2,x_3$ are distinct.
\end{proof}

For convex configurations, an axis of symmetry either passes through two  bodies or no bodies at all. We begin with the former case, that of a convex kite.

\begin{proposition}
\label{prop:convex-kite-square}
For every $A>0$, a convex kite central configuration of four equal masses is necessarily a square.
\end{proposition}

\begin{proof}
Without loss of generality, we assume
\begin{equation}
q_1=(-1,0),\quad q_2=(1,0),\quad
q_3=(0,u),\quad q_4=(0,-w),
\label{eq:convex-kite-coordinates}
\end{equation}
where $u,w>0$. Suppose first that $u\ne w$. We may assume $u>w$
and write
\[
u=m(1+k),\qquad w=m(1-k),
\qquad m>0,\quad 0<k<1.
\]
Set
\[
\begin{aligned}
a:=1+m^2(1+k)^2,\qquad b:=1+m^2(1-k)^2, \qquad e:=4m^2,
\end{aligned}
\]
and denote
\[
\psi_a:=a^{-A/2},\quad
\psi_b:=b^{-A/2},\quad
\psi_e:=e^{-A/2},\quad
\psi_4:=4^{-A/2}.
\]

The equations at $q_1$ and $q_2$ give
$\lambda=2\psi_4+\psi_a+\psi_b$.
Substituting this into the remaining equations,
we obtain
\begin{align}
2(\psi_b-\psi_a)
-k(\psi_a+\psi_b-2\psi_4)&=0,
\label{eq:convex-kite-R1}\\
k(\psi_b-\psi_a)
-2(\psi_e-\psi_4)&=0.
\label{eq:convex-kite-R2}
\end{align}
Solving these equations for $\psi_a$ and $\psi_b$
yields
\begin{equation}
\begin{aligned}
\psi_a=\psi_4+\frac{2-k}{k^2}(\psi_e-\psi_4)\quad\text{and}\quad
\psi_b=\psi_4+\frac{2+k}{k^2}(\psi_e-\psi_4).
\end{aligned}
\label{eq:convex-kite-PQ}
\end{equation}

Since $b<a$ and $A>0$, we have $\psi_b>\psi_a$.
Equation \eqref{eq:convex-kite-R2} then implies
$\psi_e>\psi_4$. Moreover, since $0<k<1$,
the coefficients of $\psi_e-\psi_4$ in
\eqref{eq:convex-kite-PQ} are both greater
than one. Hence $\psi_b>\psi_a>\psi_e$.
It follows that $b<a<e$ and 
\begin{align}
\frac{\psi_b-\psi_e}{\psi_a-\psi_e}
=\rho(k):=\frac{(1+k)(2-k)}{(1-k)(2+k)}.
\label{eq:convex-rho}
\end{align}

Since the function $s^{-A/2}$ is strictly convex on
$(0,\infty)$, the inequality for the secant slopes
with $b<a<e$ gives
\begin{equation}
\frac{\psi_b-\psi_e}{\psi_a-\psi_e}
>\frac{e-b}{e-a}.
\label{eq:convex-secant}
\end{equation}
On the other hand, direct computation gives
\begin{equation}
\frac{e-b}{e-a}-\rho(k)
=\frac{2k[1+m^2(1-k^2)]}
{(1-k)(2+k)(e-a)}>0.
\label{eq:convex-geometric-gap}
\end{equation}
This contradicts \eqref{eq:convex-rho} and
\eqref{eq:convex-secant}. Therefore $u=w$.

Now let $u=w=m$. In this case, the equations
at $q_1$ and $q_3$ give
\[
\lambda=2\psi_4+2\psi_a
=2\psi_a+2\psi_e,
\]
where $a=1+m^2$ and $e=4m^2$.
Thus $\psi_e=\psi_4$, which implies
$e=4$ and hence $m=1$. Therefore the
configuration is a square.

\end{proof}

We next consider the case in which the symmetry axis contains no body. The four bodies then form an isosceles trapezoid.

\begin{proposition}
\label{prop:trapezoid-square}
For every $A>1$, an isosceles trapezoid central configuration of four equal masses is a square.
\end{proposition}

\begin{proof}
Without loss of generality, we assume
\[
q_1=(-p,h),\quad q_2=(p,h),\quad
q_3=(q,-h),\quad q_4=(-q,-h),
\]
where $p,q,h>0$. The center of mass is the origin.
Set
\[
d^2:=(p-q)^2+4h^2,\qquad
e^2:=(p+q)^2+4h^2.
\]
The vertical component of \eqref{eq:cc} at $q_1$
gives
\begin{equation}
\lambda=2(d^{-A}+e^{-A}).
\label{eq:trapezoid-lambda}
\end{equation}
The equations for the horizontal components at $q_1$ and
$q_3$ give
\[
\begin{aligned}
\lambda p
&=2^{1-A}p^{1-A}+(p+q)e^{-A}
 +(p-q)d^{-A},\\
\lambda q
&=2^{1-A}q^{1-A}+(p+q)e^{-A}
 +(q-p)d^{-A}.
\end{aligned}
\]
Subtracting the second equation from the first and using \eqref{eq:trapezoid-lambda}, we obtain
\begin{equation}
(p-q)e^{-A}
=2^{-A}(p^{1-A}-q^{1-A}).
\label{eq:trapezoid-sign}
\end{equation}
Since $A>1$, the function $s^{1-A}$ is strictly decreasing on $(0,\infty)$. Thus, if $p\ne q$, the two sides of \eqref{eq:trapezoid-sign} have opposite signs, which is impossible. Hence $p=q$.

The configuration is therefore a rectangle. In this case, the horizontal and vertical component equations at $q_1$ give, respectively,
\[
\lambda=2(2p)^{-A}+2e^{-A},
\qquad
\lambda=2(2h)^{-A}+2e^{-A}.
\]
It follows that $p=h$, and hence the configuration is a square.
\end{proof}

We also consider central configurations with half-turn symmetry.

\begin{lemma}
\label{lem:parallelogram-square}
For every $A>0$, a noncollinear central configuration of four equal masses with half-turn symmetry is a square.
\end{lemma}

\begin{proof}
Let the center of rotation be the origin. The four positions can be written as
\[
\{x,-x,y,-y\},
\]
where $x,y\in\mathbb R^2$ are linearly independent. The central-configuration equation at $x$ gives
\[
\lambda x
=\frac{2x}{(2|x|)^A}
+\frac{x-y}{|x-y|^A}
+\frac{x+y}{|x+y|^A}.
\]
Comparing the coefficients of $y$, we obtain $|x-y|^{-A}=|x+y|^{-A}$, which implies $x\cdot y=0$.

Set $d^2:=|x|^2+|y|^2$. The equations at $x$ and $y$ then give
\[
\lambda
=2^{1-A}|x|^{-A}+2d^{-A}
=2^{1-A}|y|^{-A}+2d^{-A}.
\]
Hence $|x|=|y|$. Together with $x\cdot y=0$, this shows that the configuration is a square.
\end{proof}

Propositions~\ref{prop:convex-kite-square} and \ref{prop:trapezoid-square} settle the convex cases. Lemma~\ref{lem:parallelogram-square} also covers
configurations with half-turn symmetry. We now prove Theorem~\ref{thm:planar-classification}.

\begin{proof}[Proof of Theorem~\ref{thm:planar-classification}]
By Lemma~\ref{lem:no-three-collinear}, a noncollinear central configuration cannot contain three collinear bodies. Moreover, by the symmetry result of Albouy \cite[Propositions~1 and~2]{Albouy1995} (see also \cite[Section~3.2]{Hampton2019}), every noncollinear central configuration of four equal masses has at least one axis of reflection symmetry.

Suppose first that the configuration is concave. Its unique interior body and one of the outer bodies must lie on the symmetry axis. Hence it is a concave
kite, and Theorem~\ref{thm:concave-classification} gives exactly the two concave configurations.

The nontrivial concave kite has exactly one reflection axis. Indeed, if it had another reflection axis, both axes would pass through the unique interior body.
Their composition would give a nontrivial rotation preserving the configuration. Such a rotation must permute the three outer bodies cyclically, so they form an equilateral triangle with the fourth body at its center. This contradicts $1/3<t<2/5$.

For a convex configuration, the symmetry axis contains either two bodies or none. In the first case, the configuration is a convex kite, which is a square by Proposition~\ref{prop:convex-kite-square}. In the second case, it is an isosceles trapezoid, which is also a square by Proposition~\ref{prop:trapezoid-square}.

Finally, the generalized Moulton theorem gives a unique collinear central configuration for each ordering of the positive masses, up to translation and scaling. Since the four masses are equal, these configurations are equivalent under relabeling. This completes the proof.
\end{proof}

Theorem~\ref{thm:planar-classification} establishes the classification for $2<A\leq3$, including the range below the Newtonian exponent. For $A>3$, Hampton~\cite{Hampton2019} obtained the same four types of central configurations. Thus the classification remains unchanged as the exponent varies over $(2,\infty)$. Combining these results, we obtain the following corollary.

\begin{corollary}
\label{cor:classification-all-A}
For every $A>2$, the planar central configurations of four equal masses are exactly the four types listed in Theorem~\ref{thm:planar-classification} up to similarity and relabeling.
\end{corollary}

\begin{proof}
The case $2<A\leq3$ follows from Theorem~\ref{thm:planar-classification}, while the case $A>3$ follows from \cite{Hampton2019}.
\end{proof}

We now prove Corollary~\ref{cor-Newonian-case}, which follows from the preceding classification and the symmetry properties of the configurations.

\begin{proof}[Proof of Corollary~\ref{cor-Newonian-case}]
Let $a<-1$, so that $A=-2a>2$.
By Corollary~\ref{cor:classification-all-A}, the noncollinear central configurations are the square, the equilateral triangle with a central mass, and the unique nontrivial concave kite.

The first two configurations have more than one reflection axis, whereas the nontrivial concave kite has exactly one, as shown in the proof of Theorem~\ref{thm:planar-classification}. The collinear configuration also has more than one reflection axis by the symmetry of the equal-mass collinear solution. Hence the nontrivial concave kite is the unique central configuration with exactly one reflection axis.
\end{proof}

The classification part of the argument is now complete.  The remaining task is to determine nondegeneracy and Morse index, and to translate the geometric list into the Morse polynomial in Hampton's direct-isometry quotient.

\section{Nondegeneracy and Morse indices}
\label{sec:morse}

In this section, we study the nondegeneracy and Morse indices of the central configurations classified in Sections~\ref{sec:concave}--\ref{sec:global-classification}. Following Hampton~\cite{Hampton2019}, we consider the function
\begin{equation}
 f_A(q):=2I(q)+U_A(q),
 \quad
 I(q):=\sum_{i=1}^4|q_i|^2,
 \label{eq:Hampton-f}
\end{equation}
where $U_A$ is the homogeneous potential defined in Section~\ref{sec:reduction}. Since $U_A$ is translation invariant, every critical point of $f_A$ has its center of mass at the origin and satisfies \eqref{eq:cc} with $\lambda=4$.

In the preceding sections, the scale of the configurations was chosen without fixing $\lambda$. We first show that the Morse index obtained from the constrained Hessian in this normalization agrees with that of $f_A$ modulo rotation.

\subsection{The reduced Hessian and the scale-fixed normalization}

Let $q=(q_1,\ldots,q_4)$ be a centered central configuration satisfying
\[
\nabla U_A(q)+\lambda q=0,\quad \lambda>0.
\]
We consider the subspace
\[
\mathcal V_q:=
\left\{
E=(E_1,\ldots,E_4)\in(\mathbb R^2)^4:
\sum_{i=1}^4 E_i=0,\quad
\sum_{i=1}^4 q_i\cdot E_i=0
\right\}.
\]

For $E,G\in(\mathbb R^2)^4$, define
\[
\mathcal Q_{\lambda,A}(E,G)
:=D^2\left(U_A+\frac{\lambda}{2}I\right)(q)[E,G].
\]
Direct computation gives
\begin{equation}
\begin{aligned}
\mathcal Q_{\lambda,A}(E,G)
={}&A\sum_{i<j}r_{ij}^{-A-2}
\bigl((q_i-q_j)\cdot(E_i-E_j)\bigr)\cdot
\bigl((q_i-q_j)\cdot(G_i-G_j)\bigr)\\
&-\sum_{i<j}r_{ij}^{-A}
(E_i-E_j)\cdot(G_i-G_j)
+\lambda\sum_{i=1}^4 E_i\cdot G_i.
\end{aligned}
\label{eq:constrained-Hessian}
\end{equation}
We shall use the restriction of $\mathcal Q_{\lambda,A}$ to $\mathcal V_q$.
The following lemma relates its Morse index to the Hessian of $f_A$ at the corresponding critical point.

\begin{lemma}\label{lem:scale-index-bridge}
Let $q$ be a centered central configuration with multiplier $\lambda>0$. Set
\begin{equation}
\widehat q=sq,\quad
s=\left(\frac{\lambda}{4}\right)^{1/A}.
\label{eq:Hampton-scale}
\end{equation}
Then $\widehat q$ is a critical point of $f_A$, and
\[
n^-\bigl(D^2f_A(\widehat q)\bigr)
=
n^-\bigl(\mathcal Q_{\lambda,A}|_{\mathcal V_q}\bigr).
\]
\end{lemma}

\begin{proof}
Since $U_A$ is homogeneous of degree $2-A$, the multiplier of $sq$ is $\lambda s^{-A}$. By the choice of $s$, we have $\lambda s^{-A}=4$. Hence $\widehat q$ is a critical point of $f_A$.

For $E,G\in(\mathbb R^2)^4$, the Hessian of $f_A$ at $\widehat q$ is given by
\begin{equation}
\begin{aligned}
D^2f_A(\widehat q)[E,G]
={}&A\sum_{i<j}\widehat r_{ij}^{-A-2}
\bigl((\widehat q_i-\widehat q_j)\cdot(E_i-E_j)\bigr)
\bigl((\widehat q_i-\widehat q_j)\cdot(G_i-G_j)\bigr)\\
&-\sum_{i<j}\widehat r_{ij}^{-A}
(E_i-E_j)\cdot(G_i-G_j)
+4\sum_{i=1}^4E_i\cdot G_i,
\end{aligned}
\end{equation}
where $\widehat r_{ij}=|\widehat q_i-\widehat q_j|$. Since $\widehat q=sq$ and $\lambda=4s^A$, comparison with \eqref{eq:constrained-Hessian} gives
\begin{equation}
D^2f_A(\widehat q)[sE,sG]
=s^{2-A}\mathcal Q_{\lambda,A}(E,G).
\label{eq:Hessian-scale-factor}
\end{equation}
Since $\mathcal V_{\widehat q}=\mathcal V_q$ and $s^{2-A}>0$, the two forms have the same Morse index when restricted to $\mathcal V_q$.

We now consider the full Hessian $D^2f_A(\widehat q)$. Let
\[
\mathcal T=
\{(a,a,a,a):a\in\mathbb R^2\}.
\]
Since $\widehat q$ is centered, we have the orthogonal decomposition
\[
(\mathbb R^2)^4
=
\mathcal T\oplus\mathbb R\widehat q
\oplus\mathcal V_q.
\]

By translation invariance of $U_A$,
\[
D^2f_A(\widehat q)[(a,a,a,a),E]
=4\sum_{i=1}^4a\cdot E_i.
\]
Thus $D^2f_A(\widehat q)$ is positive definite on $\mathcal T$, and $\mathcal T$ is Hessian-orthogonal to the other two subspaces.

By homogeneity and the critical-point equation, we obtain
\begin{equation}
D^2f_A(\widehat q)[\widehat q,\widehat q]
=A(A-2)U_A(\widehat q)>0.
\label{eq:radial-positive}
\end{equation}

Moreover, for every $E\in\mathcal V_q$, homogeneity and the critical-point equation yield
\[
\begin{aligned}
D^2f_A(\widehat q)[\widehat q,E]
=(1-A)DU_A(\widehat q)[E]
  +4\sum_{i=1}^4\widehat q_i\cdot E_i
=4A\sum_{i=1}^4\widehat q_i\cdot E_i
=0.
\end{aligned}
\]
Therefore the radial direction is positive and Hessian-orthogonal to $\mathcal V_q$.

Consequently, the full Hessian and its restriction to $\mathcal V_q$ have the same number of negative eigenvalues. Together with
\eqref{eq:Hessian-scale-factor}, this completes the proof.
\end{proof}

\subsection{The nontrivial concave kite}
We now compute the Morse index of the nontrivial concave central configuration. Fix $2<A\leq3$ and consider the configuration obtained in
Theorem~\ref{thm:concave-classification}:
\[
q_1=(-1,0),\quad q_2=(1,0),\quad
q_3=(0,u),\quad q_4=(0,v),
\quad v=tu,
\]
where
\[
\frac13<t<\frac25,\quad
u^2=X(t,A)>3.
\]
Its center of mass is $c=(0,c_y)$, where
\[
c_y=\frac{u+v}{4}.
\]
We use the reflection symmetry of the configuration to simplify the Hessian. Let $R(x,y)=(-x,y)$. Since the reflection exchanges $q_1$ and $q_2$ and fixes $q_3$ and $q_4$, it acts on variations by
\[
\mathcal R(E_1,E_2,E_3,E_4)
=(RE_2,RE_1,RE_3,RE_4).
\]

Since $\mathcal R^2=I$, the space
$\mathcal V_q$ decomposes as
\[
\mathcal V_q=\mathcal V_q^+\oplus\mathcal V_q^-,
\quad
\mathcal V_q^\pm
=\{E\in\mathcal V_q:\mathcal RE=\pm E\}.
\]
Indeed, the symmetry of the configuration implies
that $\mathcal V_q$ is invariant under $\mathcal R$.
Moreover, by \eqref{eq:constrained-Hessian},
\[
\mathcal Q_{\lambda,A}(\mathcal RE,\mathcal RG)
=\mathcal Q_{\lambda,A}(E,G).
\]
Thus, for $E\in\mathcal V_q^+$ and
$G\in\mathcal V_q^-$,
\[
\mathcal Q_{\lambda,A}(E,G)
=\mathcal Q_{\lambda,A}(E,-G)=0.
\]
Hence the even and odd subspaces are
orthogonal with respect to $\mathcal Q_{\lambda,A}$.

We choose the following vectors:
\begin{align*}
o_1&=((1,0),(1,0),(0,0),(-2,0)),\\
o_2&=((0,0),(0,0),(1,0),(-1,0)),\\
e_1&=((-v,1),(v,1),(0,0),(0,-2)),\\
e_2&=((u-v,0),(-(u-v),0),(0,2),(0,-2)).
\end{align*}
The infinitesimal rotation about the center of mass is given by
\begin{equation}
\omega=
((c_y,-1),(c_y,1),
(-(u-c_y),0),(-(v-c_y),0)).
\label{eq:rotation-vector}
\end{equation}

These five vectors form a basis of the five-dimensional space $\mathcal V_q$.
More precisely,
\[
\begin{aligned}
\mathcal V_q^-
=\operatorname{span}\{o_1,o_2,\omega\},\quad
\mathcal V_q^+=\operatorname{span}\{e_1,e_2\},
\end{aligned}
\]
where $\mathcal V_q^-$ and $\mathcal V_q^+$ are the odd and even subspaces, respectively.

Since $\omega$ corresponds to rotation,
\[
\mathcal Q_{\lambda,A}(\omega,E)=0,
\quad E\in\mathcal V_q.
\]
Thus, after removing this zero direction, the Hessian has the block form
$M_-\oplus M_+$, where
\[
M_-=
\bigl(\mathcal Q_{\lambda,A}(o_i,o_j)\bigr)_{i,j=1}^2,
\quad
M_+=
\bigl(\mathcal Q_{\lambda,A}(e_i,e_j)\bigr)_{i,j=1}^2.
\]
In particular,
\[
n^-\bigl(\mathcal Q_{\lambda,A}|_{\mathcal V_q}\bigr)
=n^-(M_-)+n^-(M_+).
\]

For the calculations below, we write
\[
\begin{aligned}
P=(1+u^2)^{-A/2},\quad
Q=(1+t^2u^2)^{-A/2},\quad
\rho=[u(1-t)]^{-A},\quad
c_A=2^{1-A}.
\end{aligned}
\]
We first determine the signature of $M_-$.

\begin{lemma}
\label{lem:odd-block}
For $2<A\leq3$, let $q$ be the nontrivial concave central configuration of
Theorem~\ref{thm:concave-classification}. Then $\det M_-<0$. In particular, $M_-$ has one negative and one positive eigenvalue.
\end{lemma}

\begin{proof}
For the nontrivial concave configuration, we have $1/3<t<2/5$ and $u^2>3$. Hence $tu^2>1$. By \eqref{eq:constrained-Hessian}, the matrix $M_-=(m_{ij})_{i,j=1}^2$ has entries
\begin{align*}
m_{11}={}&6c_A+\frac{2AP}{1+u^2}
+\frac{18AQ}{1+t^2u^2}-4\rho-12Q+4P,\\
m_{12}={}&2c_A-\frac{2AP}{1+u^2}
+\frac{6AQ}{1+t^2u^2}-4\rho-4Q+4P,\\
m_{22}={}&2c_A+\frac{2AP}{1+u^2}
+\frac{2AQ}{1+t^2u^2}-4\rho.
\end{align*}
Since $F=H=0$, we have
\[
c_A=\frac{(3-t)P+(3t-1)Q}{1+t},
\quad
\rho=\frac{t(3-t)Q-(3t-1)P}{1-t^2}.
\]
Substituting these expressions into $\det M_-$ gives
\begin{equation}
\det M_-=
\frac{64A\,\mathcal B}
{(1-t^2)(1+u^2)(1+t^2u^2)},
\label{eq:odd-det-reduction}
\end{equation}
where
\begin{equation}
\mathcal B=
P^2(1+t^2u^2)-t^2Q^2(1+u^2)
+(A-1)(1-t^2)PQ.
\label{eq:Bodd}
\end{equation}

Set $z=(1+u^2)/(1+t^2u^2)$. Then $1<z<t^{-2}$ and $P/Q=z^{-A/2}$. Using $1-t^2z=(1-t^2)/(1+t^2u^2)$, we obtain
\begin{equation}
\det M_-=
\frac{64APQ}{(1+u^2)(1+t^2u^2)}
\Phi_A(t,z),
\label{eq:odd-Phi}
\end{equation}
where
\[
\Phi_A(t,z)
=(A-1)+
\frac{z^{-A/2}-t^2z^{1+A/2}}{1-t^2z}.
\]

It remains to show that $\Phi_A(t,z)<0$. By $tu^2>1$, we have $tz>1$. Let
\[
d=-\ln(t\sqrt z)>0,
\quad
\delta=\ln(tz)>0.
\]
A direct calculation gives
\begin{equation}
\frac{z^{-A/2}-t^2z^{1+A/2}}{1-t^2z}
=-\frac{\sinh((A-1)d+A\delta)}{\sinh d}.
\label{eq:odd-sinh}
\end{equation}
Since $A\geq2$ and $\delta>0$,
\[
\sinh((A-1)d+A\delta)
>\sinh((A-1)d)
\geq(A-1)\sinh d.
\]
Therefore $\Phi_A(t,z)<0$. It follows from \eqref{eq:odd-Phi} that $\det M_-<0$, which proves the lemma.
\end{proof}

We next consider the even block $M_+$. Regarding $F$ and $H$ as functions of $(u,t)$ with $x=u^2$, set
\begin{equation}
\Delta_u:=F_uH_t-F_tH_u.
\label{eq:Delta-u}
\end{equation}
Along the curve $F=0$, differentiation gives
\begin{equation}
h_t=H_t-\frac{F_t}{F_u}H_u
=\frac{\Delta_u}{F_u}.
\label{eq:ht-Delta-u}
\end{equation}
By Theorem~\ref{thm:concave-classification}, we have $h_t>0$ at the nontrivial concave configuration. Since $F_u=2uF_x<0$, it follows that
\begin{equation}
\Delta_u<0.
\label{eq:Delta-u-negative}
\end{equation}

\begin{lemma}
\label{lem:even-block}
For $2<A\leq3$, let $q$ be the nontrivial concave central configuration of Theorem~\ref{thm:concave-classification}. Then $M_+$ is positive definite.
\end{lemma}

\begin{proof}
Using \eqref{eq:constrained-Hessian} and the equations $F=H=0$, a direct calculation gives
\begin{equation}
\det M_+
=-2u\bigl[u^2(3t^2-2t+3)+8\bigr]\Delta_u.
\label{eq:even-det-identity}
\end{equation}
Since $u>0$ and $3t^2-2t+3>0$, it follows from \eqref{eq:Delta-u-negative} that
\begin{equation}
\det M_+>0.
\label{eq:even-det-positive}
\end{equation}

Moreover, $\Delta_u\ne0$, so the implicit-function theorem shows that the solution $(u,t)$ depends smoothly on $A$ near each point of the nontrivial branch. By uniqueness, these local solutions form a smooth branch for $2<A\leq3$.

We determine the sign of $(M_+)_{11}$ by considering the limit $A\downarrow2$. By Lemma~\ref{lem:upper-branch-box}, the solutions $(x(A),t(A))$ lie in $[3,4]\times[1/3,2/5]$. Proposition~\ref{prop:A2-factorization} shows that the only common zero of $F$ and $H$ in this rectangle at $A=2$ is $(3,1/3)$. Hence
\begin{equation}
(u(A),t(A))\longrightarrow
\left(\sqrt3,\frac13\right)
\quad(A\downarrow2).
\label{eq:upper-branch-endpoint-limit}
\end{equation}
Substitution into \eqref{eq:constrained-Hessian} gives
\begin{equation}
M_+\longrightarrow
\begin{pmatrix}
2&8\\
8&32
\end{pmatrix}
\quad(A\downarrow2).
\label{eq:even-endpoint-matrix}
\end{equation}
Thus $(M_+)_{11}>0$ for $A>2$ sufficiently close to $2$. By continuity, this entry remains positive along the entire branch. Indeed, if $(M_+)_{11}=0$ at some point, then
\[
\det M_+=-(M_+)_{12}^2\leq0,
\]
contrary to \eqref{eq:even-det-positive}. Therefore $(M_+)_{11}>0$ for all $2<A\leq3$. Together with $\det M_+>0$, Sylvester's criterion proves that $M_+$ is positive definite.
\end{proof}
The preceding two lemmas determine the Morse index and the kernel of the Hessian at the nontrivial concave kite.
\begin{proposition}
\label{prop:kite-index}
For every $2<A\leq3$, let $\widehat q$ be the centered and rescaled nontrivial concave kite, which is a critical point of $f_A$. Then
\[
n^-\bigl(D^2f_A(\widehat q)\bigr)=1,
\quad
\ker D^2f_A(\widehat q)
=\operatorname{span}\{\omega\},
\]
where $\omega$ is defined in \eqref{eq:rotation-vector}.
\end{proposition}

\begin{proof}
By Lemma~\ref{lem:odd-block}, $M_-$ has one negative and one positive eigenvalue, while Lemma~\ref{lem:even-block} shows that $M_+$ is positive definite. Since $\omega$ belongs to the kernel of $\mathcal Q_{\lambda,A}$,
the decomposition $0\oplus M_-\oplus M_+$ gives
\[
n^-\bigl(\mathcal Q_{\lambda,A}|_{\mathcal V_q}\bigr)=1,
\quad
\ker\bigl(\mathcal Q_{\lambda,A}|_{\mathcal V_q}\bigr)
=\operatorname{span}\{\omega\}.
\]
By the scaling identity \eqref{eq:Hessian-scale-factor} and the positive definiteness of the translation and radial directions established in
Lemma~\ref{lem:scale-index-bridge}, the same conclusions hold for $D^2f_A(\widehat q)$.
\end{proof}

\subsection{Collinear configurations}

We now consider the collinear configurations. Hampton~\cite{Hampton2019} proved that for every $A\geq2$, there is a unique collinear central configuration for each ordering of positive masses, with planar Morse index $N-2$. We also establish the nondegeneracy of the collinear configuration in the four-body case.

\begin{proposition}
\label{prop:collinear-index}
Let $A>2$ and let
\[
\widehat q=((x_1,0),(x_2,0),(x_3,0),(x_4,0)),
\quad x_1<x_2<x_3<x_4,
\]
be a collinear critical point of $f_A$ for four equal masses. Then
\[
n^-\bigl(D^2f_A(\widehat q)\bigr)=2,
\quad
\ker D^2f_A(\widehat q)
=\operatorname{span}\{\omega_{\mathrm{col}}\},
\]
where
\[
\omega_{\mathrm{col}}
=((0,x_1),(0,x_2),(0,x_3),(0,x_4)).
\]
\end{proposition}

\begin{proof}

The Morse index is two by Theorem 6 of Hampton~\cite{Hampton2019}. We prove the
assertion concerning the kernel.

Set $k_{ij}=|x_i-x_j|^{-A}$ and define
\begin{equation}
(B\xi)_i=\sum_{j\ne i}k_{ij}(\xi_j-\xi_i),
\quad \xi\in\mathbb R^4.
\label{eq:collinear-B}
\end{equation}
Then $B$ is symmetric, $B\mathbf1=0$, and
\begin{equation}
\langle\xi,B\xi\rangle
=-\sum_{i<j}k_{ij}(\xi_i-\xi_j)^2\leq0.
\label{eq:B-negative}
\end{equation}
The central-configuration equation gives
\begin{equation}
Bx=-4x,\quad x=(x_1,x_2,x_3,x_4)^T.
\label{eq:Bx}
\end{equation}

For a variation $E=(E_1,\ldots,E_4)$, write $E_i=(\xi_i,\eta_i)$, with $\xi,\eta\in\mathbb R^4$. By \eqref{eq:constrained-Hessian},
\[
\begin{aligned}
D^2f_A(\widehat q)[E,E]
=4\bigl(|\xi|^2+|\eta|^2\bigr)
+(A-1)\sum_{i<j}k_{ij}(\xi_i-\xi_j)^2-\sum_{i<j}k_{ij}(\eta_i-\eta_j)^2.
\end{aligned}
\]
Using \eqref{eq:B-negative}, we obtain
\begin{equation}
\begin{aligned}
D^2f_A(\widehat q)[E,E]
=\langle L_{\parallel}\xi,\xi\rangle+\langle L_{\perp}\eta,\eta\rangle,
\end{aligned}
\label{eq:collinear-blocks}
\end{equation}
where $L_{\parallel}=4I_4-(A-1)B$ and $L_{\perp}=4I_4+B$. Let
\[
W=\left\{\zeta\in\mathbb R^4: \sum_{i=1}^4\zeta_i=0\right\}.
\]
We define two subspaces of $\mathcal V_{\widehat q}$ by
\[
\begin{aligned}
\mathcal V_{\parallel}
&=\left\{E:E_i=(\xi_i,0),\ 1\le i\le4,\
\xi\in W\cap x^\perp\right\},\\
\mathcal V_{\perp}
&=\left\{E:E_i=(0,\eta_i),\ 1\le i\le4,\
\eta\in W\right\}.
\end{aligned}
\]
Then
\[
\mathcal V_{\widehat q}
=\mathcal V_{\parallel}\oplus\mathcal V_{\perp},
\quad
\dim\mathcal V_{\parallel}=2,\quad
\dim\mathcal V_{\perp}=3.
\]

By \eqref{eq:B-negative}, $L_{\parallel}$ is positive definite. Moreover, \eqref{eq:Bx} implies $L_{\perp}x=0$. We show that this zero eigenvalue is simple.

Since $B$ is symmetric and $B\mathbf1=0$, the space $W$ is invariant under $B$. Consider the linear map
\[
D:W\longrightarrow\mathbb R^3,\quad
D\xi=(\xi_2-\xi_1,\xi_3-\xi_2,\xi_4-\xi_3)^T.
\]
This map is an isomorphism. Let
\[
C=D(B|_W)D^{-1}.
\]
A direct calculation gives the following off-diagonal entries of $C$:
\[
\begin{aligned}
C_{12}&=(k_{23}-k_{13})+(k_{24}-k_{14}),
& C_{13}&=k_{24}-k_{14},\\
C_{21}&=k_{12}-k_{13},
& C_{23}&=k_{34}-k_{24},\\
C_{31}&=k_{13}-k_{14},
& C_{32}&=(k_{13}-k_{14})+(k_{23}-k_{24}).
\end{aligned}
\]
Since $x_1<x_2<x_3<x_4$ and $A>2$,
all six entries are strictly positive.

Choose $\sigma>4$ sufficiently large that $C+\sigma I_3$ has strictly positive entries. By \eqref{eq:Bx}, the vector
\[
d:=(x_2-x_1,x_3-x_2,x_4-x_3)^T
\]
satisfies
\[
d>0,\quad
(C+\sigma I_3)d=(\sigma-4)d.
\]
By Perron's theorem, $\sigma-4$ is a simple eigenvalue of $C+\sigma I_3$.
Consequently, $-4$ is a simple eigenvalue of $B|_W$. Since $x\in W$, we obtain
\[
\ker(L_{\perp}|_W)=\operatorname{span}\{x\}.
\]

Finally, $L_{\perp}\mathbf1=4\mathbf1$, and $\mathbb R^4=\operatorname{span}\{\mathbf1\}\oplus W$ is an invariant decomposition for $L_{\perp}$.
Therefore
\[
\ker L_{\perp}=\operatorname{span}\{x\}.
\]
Together with the positive definiteness of $L_{\parallel}$, the decomposition
\eqref{eq:collinear-blocks} yields
\[
\ker D^2f_A(\widehat q)
=\operatorname{span}\{\omega_{\mathrm{col}}\}.
\]
This completes the proof.

\end{proof}

\subsection{The high-symmetry configurations and the Morse polynomial}
It remains to consider the square and the equilateral triangle with a fourth mass at its center. Hampton~\cite{Hampton2019} computed the eigenvalues of the Hessian of $f_A$ for both configurations. For the square, the only zero eigenvalue corresponds to rotation, and all other eigenvalues are positive. For the triangle with a central mass, the eigenvalues are $0,4,4,4A$ and two double eigenvalues $\lambda_-$ and $\lambda_+$, with $\lambda_-<0<\lambda_+$ for $A>2$. Thus both configurations have no additional zero eigenvalues, and their Morse indices are zero and two, respectively.

We now count the labeled critical orbits, identifying configurations only under orientation-preserving isometries. In particular, reflections are not identified. The multiplicity of each configuration is $4!$ divided by the order of its rotational stabilizer. The resulting indices and multiplicities are
\begin{equation}
\begin{array}{l|c|c|c}
\text{configuration} & \text{stabilizer} & \text{multiplicity} & \text{index}\\
\hline
\text{square} & C_4 & 24/4=6 & 0\\
\text{nontrivial concave kite} & \{1\} & 24 & 1\\
\text{triangle plus center} & C_3 & 24/3=8 & 2\\
\text{collinear} & C_2 & 24/2=12 & 2
\end{array}
\label{eq:morse-table}
\end{equation}
Here $C_k$ denotes the cyclic group of order $k$. For the collinear configuration, uniqueness for each ordering and equality of the masses imply that $x_1=-x_4$ and $x_2=-x_3$. Hence the half-turn is a symmetry, giving the stabilizer $C_2$.

We have now determined the Morse indices and labeled multiplicities of all four configurations.
\begin{proof}[Proof of Theorem~\ref{Morse Theorem.}]
By Theorem~\ref{thm:planar-classification}, the four configurations in \eqref{eq:morse-table} exhaust all planar central configurations for $2<A\leq3$. Propositions~\ref{prop:kite-index} and~\ref{prop:collinear-index}, together with Hampton's Hessian calculations for the square and the triangle with a central mass, show that the only Hessian zero directions arise from rotation. Summing the labeled critical orbits according to their Morse indices, we obtain
\[
M_A(s)=6+24s+(8+12)s^2
      =6+24s+20s^2.
\]
\end{proof}

Hampton proved the same Morse polynomial for $A\ge3$.  Combining his theorem with Theorem~\ref{Morse Theorem.} yields the full homogeneous range.

\begin{corollary}
\label{cor:hampton-conjecture}
For every $A>2$, the Morse polynomial of the four-equal-mass planar problem is
\[
 M_A(s)=6+24s+20s^2.
\]
\end{corollary}

\begin{proof}
Theorem~\ref{Morse Theorem.} covers $2<A\le3$, and Hampton's Theorem~5 in \cite{Hampton2019} covers $A\ge3$.
\end{proof}

\appendix
\section{Exact algebraic certificates for the strict-convexity estimate}
\label{app:exact-algebra}

This appendix supplies the polynomial definitions and exact certificates used in
Lemma~\ref{prop:direct-tangent-positive}.  All calculations are over
$\mathbb Q$ on the box
\[
 1\le\gamma\le\frac32,\qquad
 \frac13\le t\le\frac25,\qquad 3\le x\le4.
\]
The two coefficient tables below are reproduced by the supplementary exact
verification scripts.

\subsection{Reconstruction of the direct-tangent polynomial}
\label{app:M-reconstruction}

Starting from \eqref{eq:F}--\eqref{eq:H}, differentiate with $c_A$ fixed and
form $\mathcal K=H_xD^2F-F_xD^2H$ as in
\eqref{eq:relative-Hessian}.  Only after differentiation substitute
\[
 c_A=\frac{(3-t)P+(3t-1)Q}{1+t},\qquad P=YQ,\qquad R=ZQ.
\]
Clearing the positive denominator gives exactly
\begin{equation}
 -\frac{K_{22}F_x-2K_{12}F_t}{4\gamma^2Q^3}
 =\frac{\mathscr M(\gamma,t,x,Y,Z)}
 {x(1-t)^2(1+t)(1+x)^2(1+t^2x)^3}.
 \label{eq:app-M-reconstruction}
\end{equation}
This identity defines $\mathscr M\in\mathbb Q[\gamma,t,x,Y,Z]$ uniquely;
its multidegree is $(2,10,4,3,1)$ and it has $348$ nonzero monomials.
The supplementary scripts reconstruct this polynomial directly from $F,H$,
so its full power-basis expansion need not be printed.

\subsection{A coefficient lower-bound principle}
\label{app:coefficient-principle}

\begin{lemma}[Shifted-power certificate]
\label{lem:shifted-power-certificate}
Let
\[
 P(z_1,\ldots,z_m)=c_0+\sum_{\alpha\ne0}c_\alpha z^\alpha
\]
be a polynomial with rational coefficients on the unit cube
$0\le z_i\le1$.  Then
\begin{equation}
 P(z)\ge c_0+\sum_{\substack{\alpha\ne0\\ c_\alpha<0}}c_\alpha.
 \label{eq:shifted-power-bound}
\end{equation}
\end{lemma}

\begin{proof}
Since $0\le z^\alpha\le1$, each negative-coefficient term is bounded below
by its coefficient, and each nonnegative-coefficient term may be discarded.
\end{proof}

For the three-variable checks below we use the affine coordinates
\begin{equation}
 \gamma=1+\frac{u}{2},
 \qquad
 t=\frac{5+v}{15},
 \qquad
 x=3+w,
 \qquad
 0\le u,v,w\le1.
 \label{eq:app-unit-cube}
\end{equation}

\subsection{The \texorpdfstring{$Z$}{Z}-monotonicity certificate}
\label{app:Z-certificate}

Write
\[
 \mathscr M_Z=aY^2+bY+d.
\]
For the leading coefficient, the exact factorization is
\[
 a=(1+t^2x)^3\mathcal A,
\]
where
\begin{align*}
 \mathcal A={}&
 \gamma^2x(2t^3-10t^2+6t+18)
 +\gamma x(-t^3+13t^2-35t+15)\\
 &+\gamma(8t^2-32t+24)+8(t-1)(x+1).
\end{align*}
Here $\mathcal A_\gamma=(t-3)\mathcal B$, where
\[
 \mathcal B=4\gamma t^2x-8\gamma tx-12\gamma x-t^2x+10tx+8t-5x-8
 \le-\frac{899}{25}.
\]
Thus $a$ is increasing in $\gamma$.  At $\gamma=1$, the polynomial
$a_x/(1+t^2x)^2$ is increasing in $x$ and is at least $4366/125$ at $x=3$.
The polynomial $a_t/(1+t^2x)^2$ is quadratic in $x$ with positive leading
coefficient; its value and $x$-derivative at $x=3$ are at least
$13459/75$ and $3637/25$, respectively.  Hence $a$ is also increasing in
$x,t$ at $\gamma=1$, and
\begin{equation}
 a\ge a\left(1,\frac13,3\right)=\frac{4096}{27}.
 \label{eq:app-a-bound}
\end{equation}

For the constant coefficient, the exact factorization is
\[
 d=t(x+1)^2\left[
 \gamma^2A_2+\gamma t(3t-1)B
 +8(1-t)(1+t^2x)C\right],
\]
where
\[
 A_2=2t^2x(t+1)(3t-1)^2(tx-t+2)\ge0,
\]
\begin{align*}
 B={}&5t^4x^2-6t^4x-10t^3x^2+16t^3x+t^2x^2+t^2x\\
 &-18tx+8t+3x-8,
\end{align*}
and
\[
 C=t^3x-3t^2+5t-1.
\]
On the box, $B_x<0$, $B\ge-14444/625$, and
$C\ge3t^3-3t^2+5t-1>0$.  Dropping $A_2\ge0$ and using
$\gamma\le3/2$ therefore gives
\[
 \frac{d}{t(x+1)^2}\ge E(t),\qquad
 E(t):=8(1-t)(1+3t^2)(3t^3-3t^2+5t-1)
 -\frac{21666}{625}t(3t-1).
\]
The derivative satisfies $E'=-2P/625$, where
\[
 P(t)=135000t^5-225000t^4+270000t^3-180000t^2+119998t-25833.
\]
Since $P(1/3)=17497/9$ and $P'>0$ on $[1/3,2/5]$,
$E(t)\ge E(2/5)=1428/625$.  Together with $t(x+1)^2\ge16/3$, this yields
\begin{equation}
 d\ge\frac{7616}{625}.
 \label{eq:app-d-bound}
\end{equation}

It remains to control the linear coefficient $b$.  Introduce
\[
 \alpha:=\gamma-1,
 \qquad
 q:=3t-1\in\left[0,\frac15\right],
 \qquad
 y:=x-3.
\]
The exact expansion has the form
\begin{equation}
 b_t=\sum_{i=0}^2\sum_{j=0}^4 C_{ij}(q)\alpha^iy^j.
 \label{eq:app-bt-expansion}
\end{equation}
If $C_{ij}(q)=\sum_k c_kq^k$, define
\[
 \delta_{ij}:=c_2-
 \sum_{\substack{k\ge3\\c_k<0}}\frac{|c_k|}{5^{k-2}}.
\]
Since $q^k\le q^2/5^{k-2}$ for $k\ge2$, one has
\[
 C_{ij}(q)\ge c_0+c_1q+\delta_{ij}q^2.
\]
Table~\ref{tab:app-bt-certificate} gives the exact values of
$c_0,c_1,\delta_{ij}$.  All are positive, so $b_t>0$.

\begin{table}[ht]
\centering
\small
\renewcommand{\arraystretch}{1.15}
\begin{tabular}{ccccc}
\hline
$i$ & $j$ & $c_0$ & $c_1$ & $\delta_{ij}$\\
\hline
2 & 4 & $128/2187$ & $2336/2187$ & $2127284/455625$ \\
2 & 3 & $6272/2187$ & $68384/2187$ & $43044734/455625$ \\
2 & 2 & $9472/243$ & $66368/243$ & $96580604/151875$ \\
2 & 1 & $14336/81$ & $76288/81$ & $89756914/50625$ \\
2 & 0 & $20480/81$ & $91136/81$ & $29839192/16875$ \\
1 & 4 & $160/729$ & $1168/729$ & $993614/151875$ \\
1 & 3 & $4192/729$ & $31024/729$ & $20259454/151875$ \\
1 & 2 & $13184/243$ & $29056/81$ & $141411152/151875$ \\
1 & 1 & $18944/81$ & $98560/81$ & $138053572/50625$ \\
1 & 0 & $10240/27$ & $38912/27$ & $48478256/16875$ \\
0 & 4 & $352/2187$ & $1168/2187$ & $853558/455625$ \\
0 & 3 & $6688/2187$ & $23536/2187$ & $5831228/151875$ \\
0 & 2 & $12416/729$ & $18368/243$ & $132632204/455625$ \\
0 & 1 & $11776/243$ & $18176/81$ & $29026118/30375$ \\
0 & 0 & $2048/27$ & $19456/81$ & $57837416/50625$ \\
\hline
\end{tabular}
\caption{Exact coefficient certificate for \eqref{eq:app-bt-expansion}.}
\label{tab:app-bt-certificate}
\end{table}

At $t=1/3$, the factorization
\[
 b=\frac{128(x+1)(x+9)^2(\gamma x-9\gamma+3)}{6561}
 \ge-\frac{54080}{729}
\]
provides the lower endpoint bound.  At $t=2/5$,
\[
 b=\frac{(x+1)(4x+25)}{1953125}\,J_{\mathrm{end}}(\gamma,x),
\]
where
\begin{align*}
 J_{\mathrm{end}}={}&8008\gamma^2x^2+91000\gamma^2x+9604\gamma x^2
 -68225\gamma x\\
 &-963000\gamma+32400x+297000.
\end{align*}
The polynomial $J_{\mathrm{end}}$ is increasing in $x$ and convex in $\gamma$.
At $x=4$, its endpoint values are $-163508$ and $-89466$, respectively.
Thus $J_{\mathrm{end}}<0$ on the box.  Using $b_t>0$ gives
\begin{equation}
 -\frac{54080}{729}\le b<0.
 \label{eq:app-b-bound}
\end{equation}
Combining \eqref{eq:app-a-bound}--\eqref{eq:app-b-bound} gives
\begin{equation}
 4ad-b^2\ge\frac{628151447552}{332150625}>0.
 \label{eq:MZ-discriminant}
\end{equation}
Thus $\mathscr M_Z=aY^2+bY+d>0$ for every real $Y$.

\subsection{The \texorpdfstring{$Y$}{Y}-layer and the terminal three-variable check}
\label{app:Y-certificate}

Set
\[
 Y_0=\frac3{20},
 \qquad
 Z_0=\frac{61}{125}.
\]
At $Z=Z_0$, write
\[
 \mathscr M=A_3Y^3+B_2Y^2+C_1Y+C_0,
 \qquad
 A_3=8x(1-t)^2(1+t^2x)^3,
\]
and define
\[
 J:=2Y_0B_2+C_1.
\]
After the affine substitution \eqref{eq:app-unit-cube}, each of
$B_2$, $J$, and
\[
 \mathscr M_0(\gamma,t,x):=
 \mathscr M(\gamma,t,x,Y_0,Z_0)
\]
has exactly $165$ nonzero monomials.  Applying
Lemma~\ref{lem:shifted-power-certificate} gives the exact certificates in
Table~\ref{tab:app-shifted-certificates}.

\begin{table}[ht]
\centering
\scriptsize
\renewcommand{\arraystretch}{1.25}
\begin{tabular}{lccc}
\hline
quantity & constant term & negative-coefficient sum & lower bound\\
\hline
$B_2$ & $1992704/30375$
& $-1441480999333/961083984375$
& $61608919000667/961083984375$\\
$J$ & $15551488/455625$
& $-9774702608047/1922167968750$
& $55833137391953/1922167968750$\\
$\mathscr M_0$ & $20834528/2278125$
& $-223878918650221/25628906250000$
& $10509521349779/25628906250000$\\
\hline
\end{tabular}
\caption{Exact shifted-power certificates on the unit cube
\eqref{eq:app-unit-cube}.}
\label{tab:app-shifted-certificates}
\end{table}

All three lower bounds are strictly positive.  In particular,
\[
 B_2>0,
 \qquad
 2Y_0B_2+C_1>0,
 \qquad
 \mathscr M(\gamma,t,x,Y_0,Z_0)>0,
\]
which are precisely the three finite algebraic inputs used in the last part
of Lemma~\ref{prop:direct-tangent-positive}.

\section*{Acknowledgments}
The authors acknowledge the use of ChatGPT (OpenAI) in exploring mathematical
ideas, developing proof strategies, and assisting with algebraic calculations
and manuscript preparation. The authors take full responsibility for the
mathematical content of this paper.

\bibliographystyle{plain}
\bibliography{references}
\end{document}